\documentclass[a4paper,final,onecolumn,10pt]{article}%
\usepackage{graphicx}
\usepackage{makeidx}
\usepackage{nomencl}
\usepackage{amsmath}
\usepackage{amssymb}
\usepackage{epsfig}
\usepackage{epstopdf}
\usepackage{subfig}
\usepackage{bm}
\usepackage{setspace}
\usepackage{algorithmic}
\usepackage{algorithm}
\usepackage{ctable}
\usepackage{booktabs}
\usepackage{multirow}
\usepackage{color}
\usepackage{colordvi}
\floatstyle{ruled} \newfloat{Algorithm}{th}{lop}
\floatstyle{ruled} \newfloat{Algorithm}{th}{lop}

\newtheorem{theorem}{Theorem}[section]

\newcommand{\qed}{\nobreak \ifvmode \relax \else
      \ifdim\lastskip<1.5em \hskip-\lastskip
      \hskip1.5em plus0em minus0.5em \fi \nobreak
      \vrule height0.75em width0.5em depth0.25em\fi}

\usepackage{amssymb,colordvi}
\usepackage{color,soul,lineno,lipsum}
\usepackage{graphicx,subfig,epsfig}
\usepackage{epstopdf}
\usepackage{authblk}
\usepackage{tabularx}
\begin{document}
\title{{A novel variational model for image registration using Gaussian curvature}}
\date{}
\author[1]{Mazlinda Ibrahim \thanks{Mazlinda.Ibrahim@liv.ac.uk}}
\author[1]{Ke Chen \thanks{K.Chen@liv.ac.uk  \ (corresponding author), \ Web http://www.liv.ac.uk/cmit}}
\author[2]{Carlos Brito-Loeza \thanks{carlos.brito@uady.mx}}
\affil[1]{Centre for Mathematical Imaging Techniques and Department of Mathematical
Sciences, The University of Liverpool, Liverpool L69 7ZL, United Kingdom.}
\affil[2]{Facultad de Matem\'{a}ticas, Universidad Aut\'{o}noma de Yucat\'{a}n, M\'{e}xico.}
\maketitle

\begin{abstract}
 Image registration is one important task in many image processing applications.
  It aims to align two or more images so that useful information can be
extracted through comparison, combination or superposition.
This is achieved by constructing an optimal transformation which ensures that the
template image becomes similar to a given reference image. Although many models exist, designing a model capable of modelling large and smooth deformation field continues to pose a challenge.
This paper proposes  a novel variational model for image registration using the Gaussian curvature as a regulariser. The model is motivated by the surface
restoration work in geometric processing [Elsey and Esedoglu, \textit{Multiscale Model. Simul.}, (2009), pp. 1549-1573]. An effective
numerical solver is provided for the model using an augmented Lagrangian method. Numerical experiments can show that the new model
 outperforms three competing models based on, respectively, a linear curvature  [Fischer and  Modersitzki, \textit{J. Math. Imaging Vis.}, (2003), pp. 81-85], the mean curvature [Chumchob, Chen and Brito, \textit{Multiscale Model. Simul.}, (2011), pp. 89-128] and the diffeomorphic demon model [Vercauteren at al., \textit{NeuroImage}, (2009), pp. 61-72] in terms of robustness and
  accuracy.
\end{abstract}

\section{Introduction}
Image registration along with image segmentation are two of the most important tasks in imaging sciences problems. Here we  focus on image registration. Much research work has been done; for an extensive overview of registration techniques see \cite{zitova2003, Pluim2003, Sotiras2013}. The methods can be classified into parametric and non-parametric image registration based on the geometric transformation. For the first category, the transformation is governed by a finite set of image features or by expanding a transformation in terms of basis functions. The second category (which is our main concern in this paper) of non-parametric image registration methods is not restricted to a certain parametrisable set. The problem is modelled as a functional minimisation via the calculus of variations. Given two images, the reference $R$ and template $T$, the functional consists of a distance measure $\mathcal{D}(T,R,\boldsymbol{u})$ and a  regularisation term $\mathcal{S}(\boldsymbol{u})$ where $\boldsymbol{u}=(u_1(\boldsymbol{x}),u_2(\boldsymbol{x}))$ is the sought displacement vector at pixel $\boldsymbol{x}\in \Omega\subset\mathbb{R}^2$. The term $\mathcal{S}(\boldsymbol{u})$ removes the ill-posedness of the minimisation problem. We use the squared $\text{L}^2$ norm of the distance measure to quantify the differences between $T$ and $R$ as follows
\begin{equation}\label{DSSD}
\mathcal{D}(T,R,\boldsymbol{u})=
\frac{1}{2}\int_{\Omega}(T(\boldsymbol{x}+
    \boldsymbol{u}(\boldsymbol{x}))-R(\boldsymbol{x}))^2 d\Omega.
\end{equation}
The distance measure in equation (\ref{DSSD}) is the sum of the squared difference (SSD) which is commonly used and optimal for mono-modal image registration with Gaussian noise.  For multi-modal image registration
where $T,R$ cannot be compared directly,
other distance measures must be used \cite{Haber2005}.
Generally, the regularisation terms are inspired by physical processes such as elasticity, diffusion and motion curvature. As such, elastic image registration was introduced in \cite{Broit1981} which assumed that objects are deformed as a rubber band.

In previous works, higher order regularisation models \cite{Fisher2004,
CCC(2011)} were found to be the most robust while the diffeomorphic demon model \cite{Vercauteren2007} offers the most physical transform in terms of (nearly) bijective mapping.
Diffusion and total variation regularisation models based on first order derivatives are less complicated to implement but
are at a disadvantage compared to higher order regularisation models
based on second order derivatives due to two reasons. First, the former methods penalise rigid displacement. They cannot properly deal with transformations with translation and rotation. Second, low order regularisation is less effective than high order one  in producing smooth  transformations  which are important in several applications including medical imaging.
%  In image denoising, higher order regularisation terms such as mean curvature and Gaussian curvature shown to outperformed low order method.
The work of \cite{Fischer2002, Fischer2003,Fisher2004} proposed a high order regularisation model named as a linear curvature, which is an approximation of the surface (mean)
curvature and the model is invariant to affine registration.
%However, in \cite{Henn2006}, the author claims that due to the Neumann boundary conditions in the model$(\text{FMC})$, the resulting transformation does not consists of an affine transformation.
This work was later refined in \cite{Henn2005,Henn2006,Henn2006(GN)} where a related approximation to the sum of the squares  of the principal curvatures was suggested and higher order boundary conditions were recommended.
Without any approximation to the mean
curvature, the works in \cite{CCC(2011),Hnung2012} developed
useful numerical algorithms for models based on the nonlinear mean curvature regularisation and observed
advantages over the linear curvature models
 for image registration; however the effect of mesh folding (bijective maps) was not considered.
The diffeomorphic demon model \cite{Vercauteren2007a} is widely used due to its property of bijective maps;
 however the bijection is not precisely imposed.
 Another useful idea  of enforcing bijection, beyond the framework we consider, is via minimising the Beltrami coefficient which measures the distortion of the quasi-conformal map of registration transforms \cite{Lui2014}.

%Read More: http://epubs.siam.org/doi/abs/10.1137/130943406?journalCode=sjisbi

In this paper we propose a high order registration model based on Gaussian curvature and hope to achieve
large and smooth deformation without mesh folding.
Although the Gaussian curvature is closely related to the mean curvature, it turns out our new model based on the
Gaussian curvature is much better.
 The motivation of the proposed model comes from two factors. Firstly, we are inspired by the work of   Elsey and Esedoglu  \cite{Elsey2009} in  geometry processing
% Denoising of surfaces or so-called surface fairing involves removing the oscillation in the surface due to the noise that appears during the construction process.
where  Gaussian curvature of the image surface is used in  a variational formulation. The authors proposed the Gaussian curvature as a natural analogue of the total variation of Rudin, Osher and Fatemi (ROF) \cite{Rof1992}
model in geometry processing. Aiming to generalise the ROF model to surface fairing where the convex shapes in 3D have similar interpretation to the monotone functions in 1D problems for the ROF model, they   showed that,  based on the Gauss Bonnet theorem, the complete analogy of the total variation regularisation for surface fairing is the energy functional of the Gaussian curvature.
A very important fact pointed out in \cite{Elsey2009} stated that the mean curvature of the surface is {\it not} a suitable choice for surface fairing because the model is not effective for preserving important features such as creases and corners on the surface (although the model is still effective for removing noise). Their claims are also supported by the work of \cite{LeeSeo2005} where the authors illustrated several advantages of   Gaussian curvature over   mean curvature and  total variation in removing  noise in 2D images. First, Gaussian curvature preserves important structures such as edges and corners. Second, only Gaussian curvature can preserve structures with low gradient. Third, the model is effective in removing noise on small scale features. Thus, we believe that Gaussian curvature is a more natural physical quantity of the surface than   mean curvature.
Here we investigate the potential
of using Gaussian curvature to construct a high order regularisation
model for non-parametric image registration of mono-modal images.

%of the Gaussian curvature as a high order regularisation
%Secondly, our motivation comes from non-parametric image registration models based on physical quantities such as elastic or curvature.

%In this framework, we consider the surface or graph in $\mathbb{R}^{d+1}$ characterized by $z=T(x,y)$ where $(x,y) \in \Omega$. Then registration problem is equivalent to the problem of finding an smooth surface $z=T(\boldsymbol{x}+\boldsymbol{u})$ that preserve the geometrical features of the original surface $z=T(x,y)$. These geometrical features are capture through the use of curvature motion of the image surface \cite{ZhuMC,ZhuMC2014}.  Similarly, we can consider an image as a collection of level sets and the curvature motion is seen as the curvature of every level line or isophote of the image \cite{Aubert2006}. Through parametrization of the level set  function $\phi$ as $\phi=z-u(x,y)$, it allows  us to write the curvature of the surface in term of the zero level set function for the displacement field $\boldsymbol{u}(\boldsymbol{x})=(u_1(x,y),u_2(x,y))$.
%\newline
%\newline

The outline of this paper is as follows.  In \S2 we review the existing models for non-parametric image registration with focus on the demon, linear and mean curvature models.
In \S3 we introduce the mathematical background of Gaussian curvature for surfaces.  In \S4 we introduce a Gaussian curvature model and a numerical solver   to  solve the resulting Euler-Lagrange equations. We show in \S5 some numerical tests including comparisons. Finally, we discuss the parameters' selection issue for our
 model in \S6  and present our conclusions in \S7.

\section{Review of Non-parametric Image Registration}
In image registration,   given two images $T$ and $R$ (which are assumed to be compactly supported and bounded operators $T,R:\Omega \subset \mathbb{R}^d \rightarrow \mathbb{R}^+$), the task is to transform $T$ to match $R$ as closely as possible.
Although we consider $d=2$ throughout this paper, with some extra modifications, this work can be extended to $d=3$.
In non-parametric image registration, the transformation is denoted by $\boldsymbol{\varphi}$ where $\boldsymbol{\varphi}$ is a vector valued function
\begin{equation}
\boldsymbol{\varphi}(\boldsymbol{x}): \mathbb{R}^2\rightarrow\mathbb{R}^2
\nonumber
\end{equation}
where  $\boldsymbol{x}=(x,y)$.
To separate the overall mapping from the displacement field, we will define
 \begin{equation}
\boldsymbol{\varphi}(\boldsymbol{x})=\boldsymbol{x}+\boldsymbol{u}(\boldsymbol{x})
\nonumber
\end{equation}
where $\boldsymbol{u}(\boldsymbol{x})$ is the displacement field. Thus, finding $\boldsymbol{u}(\boldsymbol{x})$ is equivalent to finding $\boldsymbol{\varphi}$. A non-parametric image registration model takes the form
\begin{equation}\label{minJ}
\min_{\boldsymbol{u}(\boldsymbol{x})} \mathcal{J}_{\gamma}(\boldsymbol{u}(\boldsymbol{x}))=\mathcal{D}(T,R,\boldsymbol{u}(\boldsymbol{x}))+\gamma \mathcal{S}(\boldsymbol{u}(\boldsymbol{x}))
\end{equation}
where the distance measure $\mathcal{D}$ is given as in   (\ref{DSSD})  and the choice of regulariser $\mathcal{S}(\boldsymbol{u}(\boldsymbol{x}))$ differentiates different models.
 Here $\boldsymbol{u}(\boldsymbol{x})$ is searched over  a set $\mathcal{U}$ of admissible functions that minimise $\mathcal{J}_{\gamma}$ in   (\ref{minJ}). Usually, the set $\mathcal{U}$ is a linear subspace of a Hilbert space with Euclidean scalar product.
%\begin{equation}\label{DSSD}
%\mathcal{D}(T,R,\boldsymbol{u})=\frac{1}{2}\int_{\Omega}(T(\boldsymbol{x}+\boldsymbol{u}(\boldsymbol{x}))-R(\boldsymbol{x}))^2 d\boldsymbol{x}
%\end{equation}

The force term $\boldsymbol{f}(\boldsymbol{u})$, to be used by all models, is the gradient of (\ref{DSSD}) with respect to the displacement field $\boldsymbol{u}(\boldsymbol{x})$
\begin{equation}\label{force}
\boldsymbol{f}(\boldsymbol{u})=(f_1(\boldsymbol{u}),f_2(\boldsymbol{u}))^T=(T(\boldsymbol{x}+\boldsymbol{u}(\boldsymbol{x}))-R(\boldsymbol{x}))
\nabla_{\boldsymbol{u}}T(\boldsymbol{x}+\boldsymbol{u}(\boldsymbol{x}))
\end{equation}
which is non-linear. Different regularisation terms $\mathcal{S}(\boldsymbol{u}(\boldsymbol{x}))$  will lead to different non-parametric
image registration models.
Below we list three popular models selected for tests and comparisons. % later on.

{\bf Model LC}. The first one is the linear curvature model by %Fischer and Modersitzki
  \cite{Fischer2003,Fisher2004,Mod2004,Fischer2002}, where
\begin{equation}\label{FMC}
{\mathcal S}^{\mbox{FMC}}(\boldsymbol{u})=\int_{\Omega}\Big[ (\Delta u_1)^2+(\Delta u_2)^2
             \Big] d\Omega.
\end{equation}
This term is an approximation of the surface curvature $\iota(u_l)$ through the mapping $(x,y)\rightarrow (x,y,u_{l}(x,y)),\ l=1,2$ where
 \begin{equation}\label{MCa}
\iota(u_l)=\nabla\cdot \frac{\nabla u_l }{\sqrt{|\nabla u_l|^2+1}}\
\approx\
\Delta u_l
 \end{equation}
when $|\nabla u_l|\approx 0$. The Euler Lagrange equation for (\ref{minJ}) with ${\mathcal S}^{\mbox{FMC}}$ as the regularisation term is given by
a fourth order PDE
\begin{equation}\label{ELLC}
\gamma \Delta ^2 \boldsymbol{u}+\boldsymbol{f}(\boldsymbol{u})=0
\end{equation}
with boundary conditions  $\Delta u_l=0, \nabla\Delta u_l\cdot \boldsymbol{n}=0,\ l=1,2$ and $\boldsymbol{n}$ the unit outward normal vector. The model consists of the second order derivative information of the displacement field which results in smoother deformations compared to those obtained using first order models based on elastic and diffusion energies.
It is refined in \cite{Henn2005,Henn2006,Henn2006(GN)} with nonlinear boundary conditions.
 The affine linear transformation belongs to the kernel ${\mathcal S}^{\mbox{FMC}}(\boldsymbol{u})$ which is not the case in elastic or diffusion registration.
%The solution for problem (\ref{minJ}) can be sought by either the optimise-discretise approach (i.e. the EL
%equation to be discretised by a numerical method) or the discretise-optimise approach (i.e. the discretised functional to be optimised). Either way, we obtain a non-linear system of equations, iteratively to yield the final solution.
%The second order regularisation terms are better than the first order models because they are invariant to the affine kernel.
%Henn and Witsch \cite{Henn2005,Henn2006,Henn2006(GN)} model resulted in the same Euler Lagrange equation as in  (\ref{ELLC}) but with higher order boundary conditions.
%The method can be described by the following Algorithm \ref{Alg3}.
%\begin{Algorithm}
%\caption{Linear Curvature Image Registration.}
%\begin{enumerate}
%\item Initialisation $\boldsymbol{u}(\boldsymbol{x})=\boldsymbol{0}, \gamma$.
%\item For $k=1,2,\ldots,IMAX$
%\begin{enumerate}
%\item Calculate the force term $\boldsymbol{f}$.
%\item Solve Equation (\ref{ELLC}) using iterative methods, for example semi-implicit methods or discrete cosine transform (DCT).
%\end{enumerate}
%\item End for.
%\end{enumerate}
%\label{Alg3}
%\end{Algorithm}

{\bf Model MC}.  Next is the mean curvature model \cite{CCC(2011),Hnung2012}
\begin{equation}
{\mathcal S}^{\mbox{MC}}(\boldsymbol{u})=\int_{\Omega}\Big[ \Bbbk(\iota(u_1))+\Bbbk(\iota(u_2))
     \Big] d\Omega
\nonumber
\end{equation}
where $\Bbbk(s)=\frac{1}{2}s^2$ and $\iota$ is as defined in (\ref{MCa}).
%\begin{equation}
%\kappa(u_l)=\nabla\cdot \frac{\nabla u_l}{|\nabla u_l |_{}}=\nabla\cdot \frac{\nabla u_l}{\sqrt{|\nabla u_l |^2+1}}
%\nonumber
%\end{equation}
The Euler Lagrange equation for (\ref{minJ}) with ${\mathcal S}^{\mbox{MC}}$ as the regularisation term
is given by:
\begin{equation}\label{MCEL}
\begin{aligned}
&\gamma \nabla\cdot \Big(\frac{1}{\sqrt{|\nabla u_l |^2+1}}\nabla \Bbbk'(\iota(u_l))-\frac{\nabla u_l \cdot\nabla\Bbbk'(\iota(u_l))}{(\sqrt{|\nabla u_l |^2+1})^3}\nabla u_l\Big)+f_l(\boldsymbol{u})=0,\quad l=1,2
\end{aligned}
\end{equation}
with boundary condition $\nabla u_l\cdot \boldsymbol{n}=\nabla\iota(u_l)\cdot \boldsymbol{n}=0,\ l=1,2$.
%The algorithm for solving the mean curvature model is the same used in Algorithm \ref{Alg3} but replacing Equation (\ref{ELLC}) with (\ref{MCEL}).
One can use the multigrid method to solve equation (\ref{MCEL}) as in \cite{CCC(2011)}; refer also to  \cite{Hnung2012} for multi-modality image registration work.

{\bf Model D}.  Finally Thirion \cite{Thirion1998}  introduced the so-called demon registration method where every pixel in the image acts as the demons that force a pulling and pushing action in a similar way to what Maxwell did for solving the Gibbs paradox in thermodynamics.
 The original demon registration model is a special case of diffusion registration but it has been much studied and improved
  since 1998; see \cite{Pennec1999, Mod2004, Wang2005, Lombaert2012}.
The energy functional for the basic demon method is given by
  \begin{equation}\label{demon2}
  \mathcal{S}(\boldsymbol{u})=\|R(\boldsymbol{x})-T(\boldsymbol{x}+\boldsymbol{\widetilde{u}}+\boldsymbol{u})\|^2+\frac{\sigma_i^2}{\sigma_x^2} \| \boldsymbol{u}\|^2
  \end{equation}
where $\boldsymbol{\widetilde{u}}$ is the current displacement field, $\sigma_i^2$ and $\sigma_x^2$ account  for noise on the image intensity and the spatial uncertainty respectively. Equation (\ref{demon2}) can be linearised using first order Taylor expansion,
\begin{equation}\label{demon3}
  \mathcal{J}(\boldsymbol{u})=\|R(\boldsymbol{x})-T(\boldsymbol{x}+\boldsymbol{\widetilde{u}})+J \boldsymbol{u}\|^2+\frac{\sigma_i^2}{\sigma_x^2}\| \boldsymbol{u}\|^2
\end{equation}
where $J$ is given by
\begin{equation}
J=-\frac{\nabla R+\nabla T(\boldsymbol{x}+\boldsymbol{\widetilde{u}})}{2}
\nonumber
\end{equation}
for an efficient second order minimisation. The first order condition of (\ref{demon3}) leads to the new update for $\boldsymbol{\widetilde{u}}$
\begin{equation}
\boldsymbol{u}=-\frac{R(\boldsymbol{x})-T(\boldsymbol{x}+\boldsymbol{\widetilde{u}})}{\|J\|^2+\frac{\sigma_i^2}{\sigma_x^2}}J.
\nonumber
\end{equation}
The additional use of $\boldsymbol{v}$ for $\boldsymbol{\varphi}=\exp(\boldsymbol{v})$ helps to achieve a nearly  diffeormorphic transformation
(mapping),
where $\boldsymbol{v}$ is the stationary velocity field  of the displacement field $\boldsymbol{u}$;
 see \cite{Vercauteren2009}. It should be remarked that the three main steps of the model cannot be combined into a single energy functional.

In a discrete setting, since the image domain $\Omega$ is a square,  all variational models are discretised by finite differences on a uniform grid.
Refer to \cite{CCC(2011),Mod2004}.
The vertex grid is defined by
\begin{equation}
\Omega^h=\big\{ \boldsymbol{x}_{i,j}=(x_{i},y_{j})  \  \ \big| \
\
  0\leq i\leq N_1-1,\ 0\leq j\leq N_2-1\big\}
\nonumber
\end{equation}
where we shall re-use the notation $T$ and $R$ for discrete images of size $N_1 \times N_2$.

\section{Mathematical Background  of the  Gaussian curvature}
%Differential geometry provides useful mathematical tools to develop variational  models. As such, the Monge-Ampere  equation  is a non-linear second order
%partial differential equation (PDE) that arises in the differential geometry of surfaces.
%Gaussian curvature can be considered as one   application  of the Monge-Ampere equation.

In differential geometry, the  Gaussian curvature problem  seeks to identify a hypersurface of $\mathbb{R}^{d+1}$ as a graph $z=u(\boldsymbol{x})$ over
 $\boldsymbol{x}\in \Omega \subset \mathbb{R}^d$ so that, at each point of the surface, the Gaussian curvature is prescribed.
 Let $\kappa(\boldsymbol{x})$ denote the Gaussian curvature which is a real valued function in $\Omega  \subset \mathbb{R}^d$.
 %given by $\kappa(\boldsymbol{x})$.  The Gaussian curvature at each point on the surface is given by $\kappa(\boldsymbol{x})$ and the term
 The problem is modelled by    the following equation
\begin{equation}\label{MAE}
\text{det}(D^2 u)-\kappa(\boldsymbol{x})(1+|Du|^2)^{(d+2)/2} = 0
\end{equation}
where $D$ is the first order derivative operator.
Equation  (\ref{MAE}) is one of the  Monge-Ampere equations.
 For $d=2$, we have
\begin{equation}\label{k1}
\kappa(\boldsymbol{x})\equiv -\kappa^{GC} =\frac{u_{xx}u_{yy}-u_{xy}u_{yx}}{(1+u_x^2+u_y^2)^2}.
\end{equation}

In \cite{Elsey2009}, the authors define a regularisation term using the Gaussian curvature of a closed surface based on the Gauss-Bonnet theorem.
\begin{theorem}\textbf{Gauss-Bonnet Theorem}. For a compact $C^2$ surface $\partial \Sigma$, we have
\begin{equation}
\int_{\partial \Sigma} \kappa^{GC} d\sigma=2\pi \chi
\nonumber
\end{equation}
where % $\kappa^{GC}$ is the Gaussian curvature,
$d\sigma$ is the length element to the surface and $\chi$ is the Euler characteristic of the surface.
\end{theorem}
Using this Theorem, it was shown in \cite{Elsey2009}   that  the complete analogy of the total variation regularisation for surface fairing is the energy functional of the Gaussian curvature.
The analogous term $\mathcal{S}$,  to the total variation of a function,  that appears in the ROF model \cite{Rof1992}, is given by
\begin{equation}
\mathcal{S}=\int_{\partial \Sigma } |\kappa_{}(\boldsymbol{x})| d\sigma
\nonumber
\end{equation}
where $d\sigma$ is the length element to the surface $\partial \Sigma$.

%=========
%The authors proposed the Gaussian curvature as a natural analogue of the total variation of Rudin, Osher and Fatemi (ROF) \cite{Rof1992}
%model in geometry processing. Aiming to generalise the ROF model to surface fairing where the convex shapes in 3D have similar interpretation to the monotone functions in 1D problems for the ROF model,
%=========
Gaussian curvature is one of the fundamental second order geometric properties of a surface. According to the Gauss's Theorema Egregium, Gaussian curvature  is intrinsic. For a local isometric mapping  $f:\partial \Sigma\rightarrow \partial \Sigma' $ between two surfaces, Gaussian curvature remains invariant i.e.
 if $p\in \partial \Sigma$ and $p'\in \partial \Sigma'$, then $\kappa^{GC}(p)=\kappa^{GC}(p')$ and the mapping $f$ is smooth and diffeomorphic.

We can also use a level set function to define the Gaussian curvature.
Denote by $\phi$ the zero level set of the surface generated through the mapping  $(x,y):\rightarrow(x,y,u(x,y))$. Then $\phi=u(x,y)-z$ and $\nabla \phi=(u_x,u_y,-1)^T$ where $u_x=\frac{\partial u}{\partial x}$ and $u_y=\frac{\partial u}{\partial y}$. The Gaussian curvature of the level set is given by
\begin{equation}
\kappa^{GC}=\frac{\nabla \phi H^*(\phi) \nabla \phi^T}{|\nabla \phi|^4}
\end{equation}
where $\nabla \phi=(\phi_x,\phi_y,\phi_z)^T$, $|\nabla \phi|=\sqrt{\phi_x^2+\phi_y^2+\phi_z^2}$, $H(\phi)$ is the Hessian matrix and $H^*(\phi)$ is the adjoint matrix for $H(\phi)$. We have
\begin{equation}
H(\phi)=\left[
\begin{array}{ccc}
u_{xx} & u_{xy} &0 \\
u_{yx} & u_{yy} & 0\\
0 & 0 &0
\end{array}\right],\qquad
%\nonumber
%\end{equation}
%and
%\begin{equation}
H^*(\phi)=\left[
\begin{array}{ccc}
0 & 0 &0 \\
0& 0 & 0\\
0 & 0 &u_{xx}u_{yy}-u_{yx}u_{xy}
\end{array}\right].
\nonumber
\end{equation}
Then,
\begin{equation}
\kappa^{GC}=\frac{u_{yx}u_{xy}-u_{xx}u_{yy}}{(u_x^2+u_y^2+1)^2}.
\nonumber
\end{equation}
This is why we set $\kappa^{GC}=-\kappa(\boldsymbol{x})$ in equation (\ref{k1}). We shall
use $|\kappa^{GC}|$ to measure the Gaussian curvature as in \cite{Elsey2009} for a monotonically increasing function (since the functional should be nonnegative).

\section{Image Registration based on Gaussian Curvature}
Before introducing our new image registration model, we first
illustrate some facts to support the use of
Gaussian curvature.

\subsection{Advantages of a Gaussian curvature}
{\bf The total variation and Gaussian curvature}. We use the volume based analysis introduced in \cite{LeeSeo2005} to compare two denoising
 models, respectively based on Gaussian curvature and the total variation:
\begin{equation}\label{GCDD}
\frac{\partial u}{\partial t}=\nabla \cdot \Big(\kappa \Big(\Big|\frac{u_{xy}^2-u_{xx}u_{yy}}{(u_x^2+u_y^2+1)^2}\Big|\Big)\nabla u\Big ),
\end{equation}
\begin{equation}\label{TVCC}
\frac{\partial u}{\partial t}=\nabla \cdot
 \frac{\nabla u}{|\nabla u|}.\hspace*{25mm}
\end{equation}
Consider, for each $\alpha>0$, the time change of the volume
$v_{t,\alpha}=\{(x,y,z)\ | \ 0<z<|u(x,y,t)-\alpha|\}$
which is enclosed by the surface $z=u(x,y,t)$ and
 the plane $z=\alpha$.
 Assume   $|u(x,y,t)-\alpha|=(u(x,y,t)-\alpha)s$  with $s$ either positive ($s=1$) or negative  ($s=-1$) at all points.
Denote  by $c_{t,\alpha}$ the  closed curve defined by the
level set $u(x,y,t)=\alpha$ and accordingly by $d_{t,\alpha}$
the 2D region enclosed by $c_{t,\alpha}$. The volume change in $v_{t,\alpha}$ in time is given by
\begin{equation}
V=\frac{\partial}{\partial t}\int_{v_{t,\alpha}}dzdA
=\frac{\partial}{\partial t}\int_{v_{t,\alpha}}\int_0^{|u(x,y,t)-\alpha|}dzdA =
\frac{\partial}{\partial t}\int_{d_{t,\alpha}}|u(x,y,t)-\alpha|dA
\nonumber
\end{equation}
where
$dA$ is the area element. We now consider how $V$ changes from evolving
(\ref{GCDD}) or (\ref{TVCC}).

If $u$ is the solution of equation (\ref{TVCC}), then from Gauss' theorem
\begin{equation}
\begin{aligned}
V&=\frac{\partial}{\partial t}\int_{d_{t,\alpha}}|u(x,y,t)-\alpha|dA
=s\int_{d_{t,\alpha}}\frac{\partial u}{\partial t}dA %\\ &
=s\int_{d_{t,\alpha}}\nabla \cdot \frac{\nabla u}{|\nabla u|}dA=s\int_{c_{t,\alpha}}\frac{\nabla u}{|\nabla u|}\cdot \boldsymbol{n} d\sigma
\end{aligned}
\nonumber
\end{equation}
where $d\sigma$ is the length element and $\boldsymbol{n}$ is the unit normal vector to the curve $c_{t,\alpha}$ which is   represented as
 $\boldsymbol{n}= s\frac{\nabla u}{|\nabla u|}$. Then
\begin{equation}
V= s^2\int_{c_{t,\alpha}}\frac{\nabla u}{|\nabla u|}\cdot \frac{\nabla u}{|\nabla u|} d\sigma
 =   \int_{c_{t,\alpha}}  d\sigma
=  |c_{t,\alpha}|
\nonumber
\end{equation}
   where $|c_{t,\alpha}|$ is the length of the curve $c_{t,\alpha}$.
  % Thus, $V$ depends on the curve $c_{t,\alpha}$.
   Furthermore,   the volume variation in time is
\begin{equation}
\begin{aligned}
\int_{d_{t+\delta t,\alpha}} |u(x,y,t+\delta t)-\alpha |dA&\approx \int_{d_{t,\alpha}} |u-\alpha|dA+s\delta t
     \int_{d_{t,\alpha}}\frac{\partial u}{\partial t}dA
=\int_{d_{t,\alpha}}\Big[ |u-\alpha|
    +\delta t\frac{|c_{t,\alpha}|}{|d_{t,\alpha}|}\Big]dA
\end{aligned}
\nonumber
\end{equation}
where $|d_{t,\alpha}|$ denotes the area of the region $d_{t,\alpha}$. We can see that the change in $u$ from
 $t$ to $t+\delta t$ is proportional to the ratio $\frac{|c_{t,\alpha}|}{|d_{t,\alpha}|}$. So, when this ratio is large (indicating possibly a noise presence), the total variation model reduces it and hence removes  noise.
 %or the displacement field (in the case of image registration)
 %which has  a large level set ratio.
  However, important features of $u$ which have a large level set ratio are removed also and so not preserved by the total variation model (\ref{TVCC}).

 Using similar calculations to before for the Gaussian curvature scheme (\ref{GCDD}), we have
\begin{equation}\label{GCDD1}
\begin{aligned}
V& =\frac{\partial}{\partial t}\int_{d_{t,\alpha}}|u(x,y,t)-\alpha|dA
=s\int_{d_{t,\alpha}} \nabla \cdot  \Big(\kappa \Big(\Big|\frac{u_{xy}^2-u_{xx}u_{yy}}{(u_x^2+u_y^2+1)^2}\Big|\Big)\nabla u\Big ) dA\\
&=
s\int_{c_{t,\alpha}}  \Big(\kappa \Big(\Big|\frac{u_{xy}^2-u_{xx}u_{yy}}{(u_x^2+u_y^2+1)^2}\Big|\Big)\nabla u\Big ) \cdot \boldsymbol{n} d\sigma\
 =\int_{c_{t,\alpha}}  \Big(\kappa \Big(\Big|\frac{u_{xy}^2-u_{xx}u_{yy}}{(u_x^2+u_y^2+1)^2}\Big|\Big)\Big )|\nabla u| d \sigma.\\
\end{aligned}
\end{equation}
From here, we observe that the quantity $V$ for the subdomain $v_{t,\alpha}$ is dependent on the product of the variation \textit{and} the Gaussian curvature on the level curve. The function $\kappa$ in (\ref{GCDD1}) controls and scales the speed of the volume change in contrast to the total variation scheme where $V$ depends only on the variation of the level curve.
Consider a point $p=(x_0,y_0,\alpha)$ where $\alpha=u(x_0,y_0)$. Gaussian curvature $\kappa=\kappa_1\kappa_2$ based on two principal curvatures $\kappa_1$ and $\kappa_2$ where $\kappa_1$  is the curvature of the level curve passing the point $p$ and $\kappa_2$ is the curvature of the path which passes the point and $\kappa_2$ is orthogonal to the level curve. If the Gaussian curvature on one level curve is zero then there is no change in $V$ regardless of variation  on the level curves. In contrast, with the total variation, if there is a variation in the level curve, then there is a change in $V$. Based on this observation, we believe that the Gaussian curvature model is better than the total variation model for preserving features on surfaces.
%%%%compare mc and gc

{\bf The mean curvature and Gaussian curvature}. The mean curvature (MC)  $\iota=(\kappa_1+\kappa_2)/2$ is also widely used. Next, we show that, though closely related,
 Gaussian curvature (GC) is better than mean curvature for surfaces in three ways.

 First, Gaussian curvature is invariant under rigid and isometric transformations. In contrast, mean curvature is invariant under rigid transformations but not under isometric transformations. Rigid transformations preserve distance between two points while isometric transformations preserve length along surfaces and preserve angles between curves on surfaces. To illustrate invariance, consider a surface
\begin{equation}
z_1(x,y)=ax^2+by^2,
\nonumber
\end{equation}
whose
Gaussian curvature and   mean curvature are respectively
\begin{equation}
\kappa=\frac{0-(2a)(2b)}{(1+4a^2x^2+4b^2y^2)^2}, \qquad
\iota=\frac{(1+4b^2y^2)(2a)+(1+4a^2x^2)(2b)}{(1+4a^2x^2+ab^2y^2)^{3/2}}.
\nonumber
\end{equation}
If we flip the surface upside down (isometric transformation) where $z_{1}'(x,y)=-ax^2-by^2$, we will have the same value for the Gaussian curvature and a different value for the mean curvature. Thus, Gaussian curvature is invariant under isometric transformation.

Second, Gaussian curvature can be used to localise the
tip of a surface better than mean curvature. Consider
\begin{equation}
z_2(x,y)=-\frac{1}{2}(x^2+y^2).
\nonumber
\end{equation}
When we compute the mean and Gaussian curvature for the surface,
 we see that the value given by the Gaussian curvature is sharper than that of the mean curvature. The highest point of the Gaussian curvature is better distinguished from its neighbourhood compared to the highest point of the mean curvature.

Third, Gaussian curvature can   locate   saddle points better than mean curvature.  Take
\begin{equation}
z_3(x,y)=-\frac{1}{2}(x^2-y^2)
\nonumber
\end{equation}
as one example.
We can again compute its mean and Gaussian curvatures analytically.
The mean curvature for this surface appears complex where the largest value is not at the saddle point and the saddle point cannot be easily located. However, Gaussian curvature gets its highest value at the saddle point and is therefore able to accurately identify the saddle point within its neighbourhood.
%$z_2(x,y)$ and $z_3(x,y)$, we can argue that Gaussian curvature can be used to locate tip and saddle points better than mean curvature.

%In addition to   these three examples and observations,
%a very important fact pointed out in \cite{Elsey2009} stated that the mean curvature of the surface is {\it not} a suitable choice for surface fairing because the model is not effective for preserving important features such as creases and corners on the surface (although the model is still effective for removing noise). This claim is also supported by the work of \cite{LeeSeo2005} where several advantages of the Gaussian curvature over the mean curvature and the total variation are illustrated  in removing  noise in 2D images.
%
%Therefore, we might conjecture that  Gaussian curvature may outperform
% the mean curvature in image registration.
%To our knowledge there exists no previous work on this topic.
In addition to these three examples and observations, a very important fact pointed out in \cite{Elsey2009} states that the mean curvature of the surface is not a suitable choice for surface fairing because the model is not effective for preserving important features such as creases and corners on the surface (although the model is effective for removing noise). This is true when we   refer  to surface fairing (surface denoising) but not necessarily true for 2D image denoising. From the recent works done in image denoising \cite{Elsey2009,LeeSeo2005}, we observe  several advantages of Gaussian curvature over total variation and mean curvature. Therefore, we   conjecture that Gaussian curvature may outperform existing models in image registration. To our knowledge there exists no previous work on this topic.
\subsection{The proposed registration model}
Now we return to the problem of how to align or register two image functions
$T(\boldsymbol{x}), R(\boldsymbol{x})$.
Let the desired and unknown displacement fields between $T$ and $R$
be the surface map $(x,y):\rightarrow(x,y,u_l(x,y))$ where $ \ l=1,2$ and with $\boldsymbol{u}=(u_1, u_2)$.
 We propose our Gaussian curvature based image registration model as
\begin{equation}\label{gc}
\begin{aligned}
& \min_{\boldsymbol{u}\in C^2(\Omega)} \mathcal{J}_{\gamma}(\boldsymbol{u}(\boldsymbol{x}))=\frac{1}{2} \int_{\Omega} \left(T(\boldsymbol{x}+\boldsymbol{u})-R(\boldsymbol{x})\right)^2 d\Omega +\gamma \mathcal{S}^{\text{GC}}(\boldsymbol{u}(\boldsymbol{x}))\\
\end{aligned}
\end{equation}
where
 \begin{equation}
\begin{aligned}
\mathcal{S}^{\text{GC}}(\boldsymbol{u}(\boldsymbol{x}))=\sum_{l=1}^2 {\mathcal S}^{GC}(u_l),\quad {\mathcal S}^{GC}(u_l)=\int_{\Omega} \left| \frac{u_{l,xy}u_{l,yx}-u_{l,xx}u_{l,yy}}{(u_{l,x}^2+u_{l,y}^2+1)^2}\right|d\Omega.
 \end{aligned}
\nonumber
\end{equation}
The above model (\ref{gc}) leads to  two Euler Lagrange equations:
\begin{equation}\label{EL}\left\{
\begin{aligned}
&\gamma\nabla \cdot \left(\frac{4|u_{1,xy}u_{1,yx}-u_{1,xx}u_{1,yy}|}{\boldsymbol{\mathcal{N}}_1^3}\nabla u_{1}\right)+\gamma\nabla \cdot \boldsymbol{B}_{1,1}+\gamma\nabla \cdot \boldsymbol{B}_{1,2}+ f_{1}=0\\
&\gamma\nabla \cdot \left(\frac{4|u_{2,xy}u_{2,yx}-u_{2,xx}u_{2,yy}|}{\boldsymbol{\mathcal{N}}_2^3}\nabla u_{2}\right)+\gamma\nabla \cdot \boldsymbol{B}_{2,1}+\gamma\nabla \cdot \boldsymbol{B}_{2,2}+f_{2}=0\\
\end{aligned}\right.
\end{equation}
where
\begin{equation}
\begin{aligned}
&\boldsymbol{\mathcal{N}}_l=u_{l,x}^2+u_{l,y}^2+1, \ \ \boldsymbol{B}_{l,1}=\left(\left(-\frac{\mathcal{S}_lu_{l,yy}}{\boldsymbol{\mathcal{N}}_{l}}\right)_{x},\left(\frac{\mathcal{S}_lu_{l,xy}}{\boldsymbol{\mathcal{N}}_{l}}\right)_{x}\right)\\
&\boldsymbol{B}_{l,2}=\left(\left(\frac{\mathcal{S}_lu_{l,yx}}{\boldsymbol{\mathcal{N}}_{l}}\right)_{y},\left(-\frac{\mathcal{S}_lu_{l,xx}}{\boldsymbol{\mathcal{N}}_{l}}\right)_{y}\right),
\ \
\mathcal{S}_l=\text{sign}(u_{l,xy}u_{l,yx}-u_{l,xx}u_{l,yy})\\
& \boldsymbol{f}=(f_1,f_2)^T=(T(\boldsymbol{x}+\boldsymbol{u})-R(\boldsymbol{x}))\nabla_{\boldsymbol{u}}T(\boldsymbol{x}+\boldsymbol{u}),
\ l=1,2
\end{aligned}
\nonumber
\end{equation}
and boundary conditions $\nabla u_l \cdot \boldsymbol{n}=0,\ l=1,2$, where $\boldsymbol{n}$ is the normal vector at the boundary $\partial \Omega$.
%$\boldsymbol{f}$ is the first variation of the similarity measure $D(T,R,\boldsymbol{u})=\int_{\Omega} \left(T(\boldsymbol{x}+\boldsymbol{u})-R(\boldsymbol{x})\right)^2 d\Omega$.
Derivation of the resulting Euler-Lagrange equations for this model can be found in Appendix A. The equations are fourth order and nonlinear with anisotropic diffusion.
 Gradient descent type methods are not feasible and
 one way to solve them is by a geometric multigrid as in \cite{CCC(2011)}. Here, instead,  we present a fast and efficient way to  solve the model (\ref{gc}) on a  unilevel grid.

\subsection{Augmented Lagrangian Method}
The augmented Lagrangian method (ALM) is often used for solving constraint minimisation problems by replacing the original problem with an unconstrained problem. The method is similar to the penalty method where the constraints are incorporated in the objective functional and the problem is solved using alternating minimisation of each of the sub-problems. However, in ALM, there are additional terms in the objective functional known as Lagrange multiplier terms arising when incorporating the constraints.
Similar works on the augmented Lagrangian method in image restoration  can be found in \cite{ZhuMC,ZhuMC2014}. %where the method is used to solve the mean curvature model for image denoising.

%The Euler-Lagrange Equations (\ref{EL}) are fourth order and therefore difficult to solve compared to the lower order PDEs.
To proceed, we introduce two new dual variables $\boldsymbol{q}_1$ and $\boldsymbol{q}_2$ where $\boldsymbol{q}_1=\nabla u_1(\boldsymbol{x})$ and $\boldsymbol{q}_2=\nabla u_2(\boldsymbol{x})$.
Consequently we obtain
 a system of second order PDEs which are more amendable to effective solution.
 %Denote by $W$  the Euclidean space $\mathbb{R}^{m\times n}$ of matrices of size $m \times n$ and $Q=W\times W$ where $\nabla u_l \in Q, l=1,2$. The space $Q$ can be considered as the space for the gradient of the displacement field $\boldsymbol{u}(\boldsymbol{x})=(u_1(\boldsymbol{x}),u_2(\boldsymbol{x}))$.

We obtain the following refined model for Gaussian curvature image registration
\begin{equation}
\begin{aligned}
\min_{u_1,u_2, \boldsymbol{q}_1,\boldsymbol{q}_2}  \mathcal{J}(u_1,u_2,\boldsymbol{q}_1,\boldsymbol{q}_2)&
=\mathcal{D}(T,R,\boldsymbol{u}(\boldsymbol{x}))+\gamma {\mathcal S}^{GC}(\boldsymbol{q}_1)+\gamma {\mathcal S}^{GC}(\boldsymbol{q}_2)\\
& \text{s.t} \hspace{0.1cm} \boldsymbol{q}_1=\nabla u_1(\boldsymbol{x}),\boldsymbol{q}_2=\nabla u_2(\boldsymbol{x})
\end{aligned}
\nonumber
\end{equation}
%where $D^{SSD}(T,R,u_1,u_2)=\ \int_{\Omega} \left(T(\boldsymbol{x}+\boldsymbol{u})-R(\boldsymbol{x})\right)^2 d\Omega$.
and further reformulate $\mathcal{J}(u_1,u_2,\boldsymbol{q}_1,\boldsymbol{q}_2)$ to get the augmented Lagrangian functional
 \begin{equation}\label{J_2}
 \begin{aligned}
 \mathcal{L}^{GC}(u_1,u_2,\boldsymbol{q}_1,\boldsymbol{q}_2; \boldsymbol{\mu}_1,\boldsymbol{\mu}_2)=&
 \frac{1}{2}\|T(\boldsymbol{x}+
 \boldsymbol{u}(\boldsymbol{x}))-R(\boldsymbol{x})\|^2_2+\gamma {\mathcal S}^{GC}(\boldsymbol{q}_1)+\gamma {\mathcal S}^{GC}(\boldsymbol{q}_2)\\
 &+\langle\boldsymbol{\mu}_1,\boldsymbol{q}_1-\nabla u_1\rangle+\langle\boldsymbol{\mu}_2,\boldsymbol{q}_2-\nabla u_2\rangle\\
 &+\frac{r}{2}\|\boldsymbol{q}_1-\nabla u_1\|^{2}_2+\frac{r}{2}\|\boldsymbol{q}_2-\nabla u_2\|^{2}_2
 \end{aligned}
 \end{equation}
where $\boldsymbol{\mu}_1,\boldsymbol{\mu}_2$ are the Lagrange multipliers,
the inner products are defined via the usual
integration in $\Omega$ and $r$ is a positive constant.
 %The  saddle point problem for  the augmented Lagrangian method for the Gaussian curvature image registration is
% \begin{equation}
% \begin{aligned}
% &\text{Find  } (u_1,u_2,\boldsymbol{q}_1,\boldsymbol{q}_2,\mu_1,\mu_2) \in W \times Q \times Q\\
% & \text{s.t}  \mathcal{L}^{GC}(u_1,u_2,\boldsymbol{q}_1,\boldsymbol{q}_2; \boldsymbol{\mu}_1,\boldsymbol{\mu}_2)
% \end{aligned}
% \end{equation}
 We use an alternating minimisation procedure to find the optimal values of $u_1,u_2,\boldsymbol{q}_1,\boldsymbol{q}_2$ and $\boldsymbol{\mu}_1,\boldsymbol{\mu}_2$ where the process involves only two main steps.

 {\bf Step 1}. For the first step we need to update $\boldsymbol{q}_1,\boldsymbol{q}_2$ for any given $u_1,u_2,\boldsymbol{\mu}_1,\boldsymbol{\mu}_2$. The objective functional is given by
\begin{equation}
\begin{aligned}
\min_{\boldsymbol{q}_1,\boldsymbol{q}_2} &\gamma {\mathcal S}^{GC} (\boldsymbol{q}_1)+\gamma {\mathcal S}^{GC} (\boldsymbol{q}_2)+\langle\boldsymbol{\mu}_1,\boldsymbol{q}_1\rangle
+\langle\boldsymbol{\mu}_2,\boldsymbol{q}_2\rangle
+\frac{r}{2}\|\boldsymbol{q}_1-\nabla u_1\|^2+\frac{r}{2}\|\boldsymbol{q}_2-\nabla u_2\|^2.
\end{aligned} \nonumber
\end{equation}
This sub-problem can be solved using the following Euler Lagrange equations:
\begin{equation}\label{EL2a}\left\{
\begin{aligned}
&  -\gamma\Big (\Big (\frac{(-q_{1,1})_y}{\Gamma_1^2}\Big)_x+ \Big (\frac{-(q_{1,1})_x}{\Gamma_1^2}\Big)_y\Big)-\gamma\frac{4S_1 D_1 q_{1,2}}{\Gamma_1^3}+\mu_{1,2}+r(q_{1,2}-(u_{1})_y)=0,\\
&  -\gamma\Big (\Big (\frac{(q_{1,2})_y}{\Gamma_1^2}\Big)_x+ \Big (\frac{-(q_{1,2})_x}{\Gamma_1^2}\Big)_y\Big)-\gamma\frac{4S_1 D_1 q_{1,1}}{\Gamma_1^3}+\mu_{1,1}+r(q_{1,1}-(u_{1})_x)=0
\end{aligned}\right.
\end{equation}
where $D_1=\text{det}(\nabla \boldsymbol{q}_1)=(q_{1,1})_x (q_{1,2})_y-(q_{1,1})_y(q_{1,2})_x$, $\Gamma_1=1+u_{1,x}^2+u_{1,y}^2$ and $S_1=\text{sign}\Big (\frac{D_1}{(\|\nabla u_1\|^2+1)^2}   \Big)$.
We have a closed form solution for this step, if solving alternatingly,  where
\begin{equation}\nonumber %\label{closedform}
\begin{aligned}
q_{1,1}&=\frac{\Gamma_1^3 \Big(-\gamma\Big (\Big (\frac{(q_{1,2})_y}{\Gamma_1^2}\Big)_x+ \Big (\frac{-(q_{1,2})_x}{\Gamma_1^2}\Big)_y\Big)+ \mu_{1,1}-r(u_1)_x\Big)}{-r \Gamma_1^3+ \gamma4S_1D_1 },\\
q_{1,2}&=\frac{\Gamma_1^3 \Big(-\gamma\Big (\Big (\frac{(q_{1,1})_y}{\Gamma_1^2}\Big)_x+ \Big (\frac{-(q_{1,1})_x}{\Gamma_1^2}\Big)_y\Big)+ \mu_{1,2}-r(u_1)_y\Big)}{-r \Gamma_1^3+ \gamma4S_1D_1 }.
\end{aligned}
\end{equation}
Similarly, we solve $q_{2,1}, q_{2,2}$ from
\begin{equation}\label{EL2b}\left\{
\begin{aligned}
&  -\gamma\Big (\Big (\frac{(-q_{2,1})_y}{\Gamma_2^2}\Big)_x+ \Big (\frac{-(q_{2,1})_x}{\Gamma_2^2}\Big)_y\Big)-\gamma\frac{4S_2 D_2 q_{2,2}}{\Gamma_2^3}+\mu_{2,1}+r(q_{2,2}-(u_{2})_y)=0,\\
&  -\gamma\Big (\Big (\frac{(q_{2,2})_y}{\Gamma_2^2}\Big)_x+ \Big (\frac{-(q_{2,2})_x}{\Gamma_2^2}\Big)_y\Big)-\gamma\frac{4S_2 D_2 q_{2,1}}{\Gamma_2^3}+\mu_{2,1}+r(q_{2,1}-(u_{2})_x)=0\\
\end{aligned}\right.
\end{equation}
where $D_2=\text{det}(\nabla \boldsymbol{q}_2)=(q_{2,1})_x (q_{2,2})_y-(q_{2,1})_y(q_{2,2})_x$, $\Gamma_2=1+u_{2,x}^2+u_{2,y}^2$ and $S_2=\text{sign}\Big (\frac{D_2}{(\|\nabla u_2\|^2+1)^2}   \Big)$.

 {\bf Step 2}.
For the second step we need to update $u_1,u_2$ for any given $\boldsymbol{q}_1,\boldsymbol{q}_2$ and $\boldsymbol{\mu}_1,\boldsymbol{\mu}_2$ with the following functional
\begin{equation}
\begin{aligned}
\min_{u_1,u_2} \frac{1}{2} \|T(\boldsymbol{x}+\boldsymbol{u})-R(\boldsymbol{x})\|_2^2-
\langle\boldsymbol{\mu}_1, \nabla u_1\rangle-
\langle\boldsymbol{\mu}_2, \nabla u_2\rangle
+\frac{r}{2}\|\boldsymbol{q}_1-\nabla u_1\|^2+\frac{r}{2}\|\boldsymbol{q}_2-\nabla u_2\|^2.
\end{aligned}
\nonumber
\end{equation}
Thus, we have the following Euler Lagrange equations:
\begin{equation}\label{EL3}\left\{
\begin{aligned}
&-r \Delta u_1+f_1+\nabla \cdot \boldsymbol{\mu}_1+r \nabla \cdot \boldsymbol{q}_1=0\\
&-r \Delta u_2+ f_2+\nabla \cdot \boldsymbol{\mu}_2+r \nabla \cdot \boldsymbol{q}_2=0\\
\end{aligned}\right.
\end{equation}
with Neumann boundary conditions $\nabla u_l \cdot \boldsymbol{n}=0, l=1,2$. To solve equation (\ref{EL3}), first, we linearise the force term $\boldsymbol{f}$ using Taylor expansion
\begin{equation}
\begin{aligned}
f_l(u_1^{(k+1)},u_2^{(k+1)})&=f_l(u_1^{(k)},u_2^{(k)})+\partial_{u_1}f_l(u_1^{(k)},u_2^{(k)})\delta u_l^{(k)}+\partial_{u_2}f_1(u_1^{(k)},u_2^{(k)})\delta u_2^{(k)}+\ldots\\
&\approx f_l(u_1^{(k)},u_2^{(k)})+\sigma^{(k)}_{l,1}\delta u_1^{(k)}+\sigma_{l,2}\delta u_2^{(k)}
\end{aligned}
\nonumber
\end{equation}
where
\begin{equation}
\begin{aligned}
\sigma^{(k)}_{l,1}&=\partial_{u_1}f_l(u_1^{(k)},u_2^{(k)}),\sigma_{l,2}=\partial_{u_2}f_l(u_1^{(k)},u_2^{(k)}),\delta u_1^{(k)}=u_1^{(k+1)}-u_1^{(k)},\delta u_2^{(k)}=u_2^{(k+1)}-u_2^{(k)}.\\
\end{aligned}
\nonumber
\end{equation}
Second, we approximate $\sigma^{(k)}_{l,1}$ and $\sigma^{(k)}_{l,2}$ with
\begin{equation}
\begin{aligned}
\sigma^{(k)}_{l,1}&=\Big(\partial u_l T(\boldsymbol{x}+\boldsymbol{u}^{(k)})\Big)\Big(\partial u_1 T(\boldsymbol{x}+\boldsymbol{u}^{(k)})\Big)\\
\sigma^{(k)}_{l,2}&=\Big(\partial u_l T(\boldsymbol{x}+\boldsymbol{u}^{(k)})\Big)\Big(\partial u_2 T(\boldsymbol{x}+\boldsymbol{u}^{(k)})\Big).\\
\end{aligned}
\nonumber
\end{equation}
The discrete version of equation (\ref{EL3}) is as follows
\begin{equation}\label{S1}
\begin{aligned}
\boldsymbol{N}^h(\boldsymbol{u}^{h,(k)})\boldsymbol{u}^{h,(k+1)}=\boldsymbol{B}^h(\boldsymbol{u}^{h,(k)})
\end{aligned}
\end{equation}
where
\begin{equation}
\begin{aligned}
\boldsymbol{N}^h(\boldsymbol{u}^{(k)})&=\Big[
\begin{array}{cc}
-r \mathcal{L}+\sigma_{11}^h(\boldsymbol{u}^{h,(k)})& \sigma_{12}^h(\boldsymbol{u}^{h,(k)})\\
 \sigma_{21}^h(\boldsymbol{u}^{h,(k)})&-r \mathcal{L}+\sigma_{22}^h(\boldsymbol{u}^{h,(k)})\\
\end{array}\Big],\\
\boldsymbol{B}^h(\boldsymbol{u}^{(k)})&=\Big[
\begin{array}{c}
-G^h_1+f^h_1(u_1^{(k)},u_2^{(k)})+\sigma_{11}^h(\boldsymbol{u}^{(k)})u_1^{h,(k)}+\sigma_{12}^h(\boldsymbol{u}^{h,(k)})u_2^{h,(k)}\\
-G^h_2+f^h_2(u_1^{(k)},u_2^{(k)})+\sigma_{21}^h(\boldsymbol{u}^{(k)})u_1^{h,(k)}+\sigma_{22}^h(\boldsymbol{u}^{h,(k)})u_2^{h,(k)}\\
\end{array}\Big],
\end{aligned}
\nonumber
\end{equation}
$\mathcal{L}$ is the discrete version of the Laplace operator $\Delta$ and $G^h_l$ is the discrete version of
\begin{equation}
\nabla \cdot \boldsymbol{\mu}_l+r\nabla \cdot \boldsymbol{q}_l,\ l=1,2.
\nonumber
\end{equation}
Third, we solve the system of equation (\ref{S1}) using a weighted pointwise Gauss Siedel method
\begin{equation}
\boldsymbol{u}^{h,(k+1)}=(1-\omega)\boldsymbol{u}^{h,(k)}+\omega \Big(\boldsymbol{N}^h(\boldsymbol{u}^{(k)})\Big)^{-1}\boldsymbol{B}^h(\boldsymbol{u}^{(k)})
\nonumber
\end{equation}
where $\omega \in(0,2)$ and we choose $\omega=0.9725$.

The iterative algorithm to solve (\ref{J_2}) is now summarised as follows.
\begin{Algorithm}
\caption{Augmented Lagrangian method for the Gaussian Curvature Image Registration.}
\begin{enumerate}
\item Initialise $\boldsymbol{\mu}_1=\boldsymbol{\mu}_2=\boldsymbol{0}, \boldsymbol{u}(\boldsymbol{x})=0, \gamma, r$.
\item For $k=0,1,...,IMAX$
\begin{enumerate}
\item Step 1: Solve
 (\ref{EL2a}-\ref{EL2b}) for $(\boldsymbol{q}_1^{(k+1)},\boldsymbol{q}_2^{(k+1)})$ with
      $(u_1, u_2)=(u_1^{(k)},u_2^{(k)})$.

\item Step 2: Solve (\ref{EL3}) for $(u_1^{(k+1)},u_2^{(k+1)})$ with $(\boldsymbol{q}_1,\boldsymbol{q}_2)=(\boldsymbol{q}_1^{(k+1)},\boldsymbol{q}_2^{(k+1)})$.

\item Step 3: Update Lagrange multipliers.
\newline
 $\boldsymbol{\mu}_1^{(k+1)}=\boldsymbol{\mu}_1^{(k)}+r(\boldsymbol{q}_1^{(k+1)}-\nabla u_1^{(k+1)})$, $\boldsymbol{\mu}_2^{(k+1)}=\boldsymbol{\mu}_2^{(k)}+r(\boldsymbol{q}_2^{(k+1)}-\nabla u_2^{(k+1)})$
\end{enumerate}
\item End for.
\end{enumerate}\label{Alg1}
\end{Algorithm}

\section{Numerical Results}
We use two numerical experiments to examine the efficiency and robustness of the Algorithm \ref{Alg1} on a variety of deformations. To judge the quality of the alignment we calculate the relative reduction of the similarity measure
\begin{equation}
\varepsilon=\frac{\mathcal{D}(T,R,\boldsymbol{u})}{ \mathcal{D}(T,R)}
\nonumber
\end{equation}
and the minimum value of the determinant of the Jacobian matrix $J$ of the transformation, denoted as ${\cal F}$
\begin{equation}
\begin{aligned}
J&=\left[
\begin{array}{cc}
1+u_{1,x} & u_{1,y}\\
u_{2,x} & 1+u_{2,y}
\end{array}\right],
\hspace{0.1cm}{\cal F}=\text{min} \left(\det(J)\right).\\
%MFN&=\text{no. of pixels with  } \text{det}(J)<0\\
\end{aligned}
\end{equation}
We can observe that when ${\cal F}>0$, the deformed grid is free from folding and cracking.

All experiments were run on a single level. Experimentally, we found that $r\in[0.02,2]$ works well for several types of images. As for the stopping criterion, we use $tol=0.001$ for the residual of the Euler-Lagrange equations  (\ref{EL2a})-(\ref{EL3}) and the maximum number of iterations is 30. Experiments were carried out using Matlab R2014b with Intel(R) core (TM) i7-2600 processor and 16G RAM.
\subsection{Test 1: A Pair of Smooth X-ray Images}
\begin{figure}[h!]
\centering
\def\tabularxcolumn#1{m{#1}}
\begin{tabular}{ccc}%\hspace*{-6mm}
 \subfloat[$T$]{\includegraphics[height=4cm,width=4.5cm]{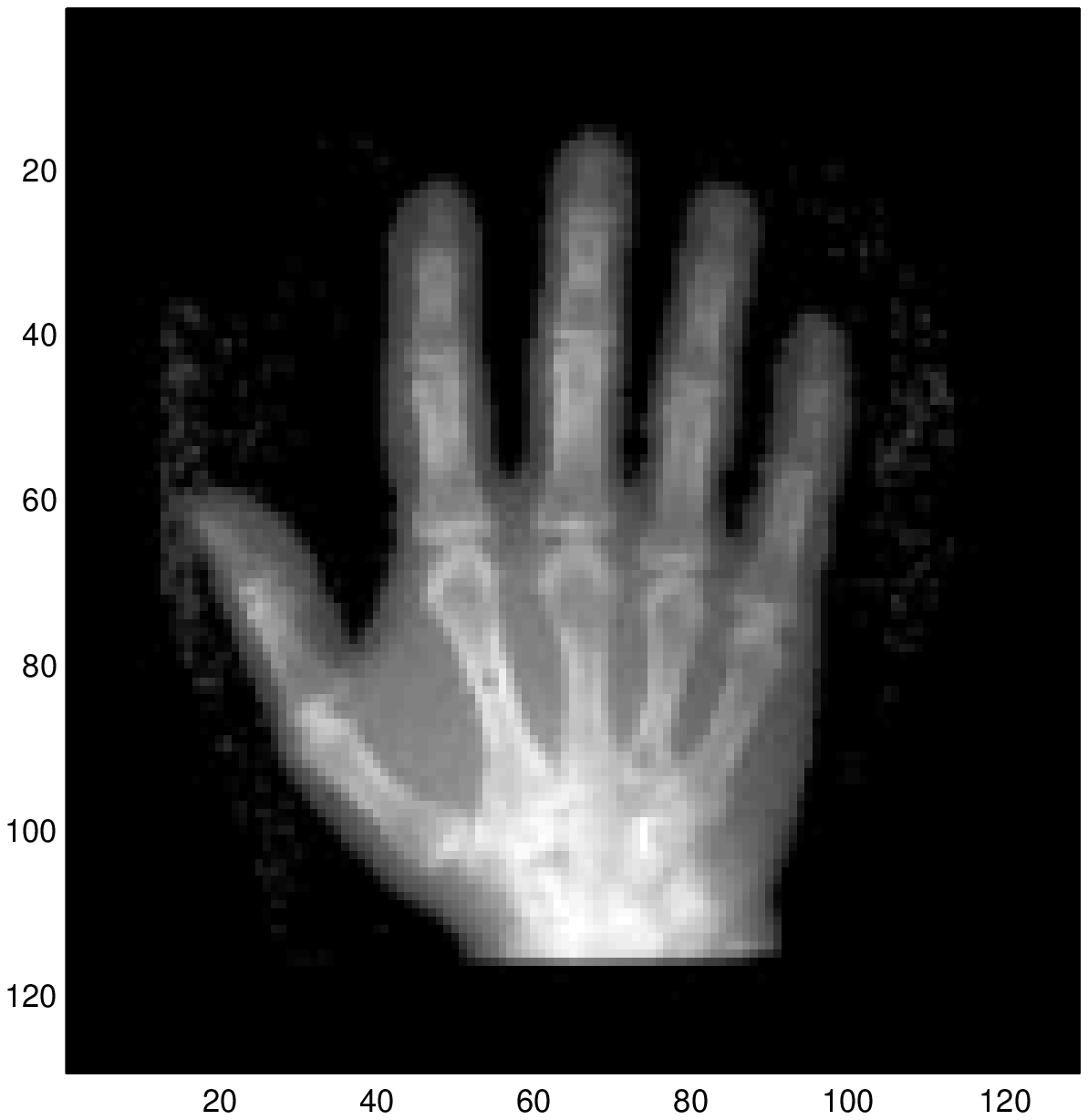}}
   & \subfloat[$R$]{\includegraphics[height=4cm,width=4.5cm]{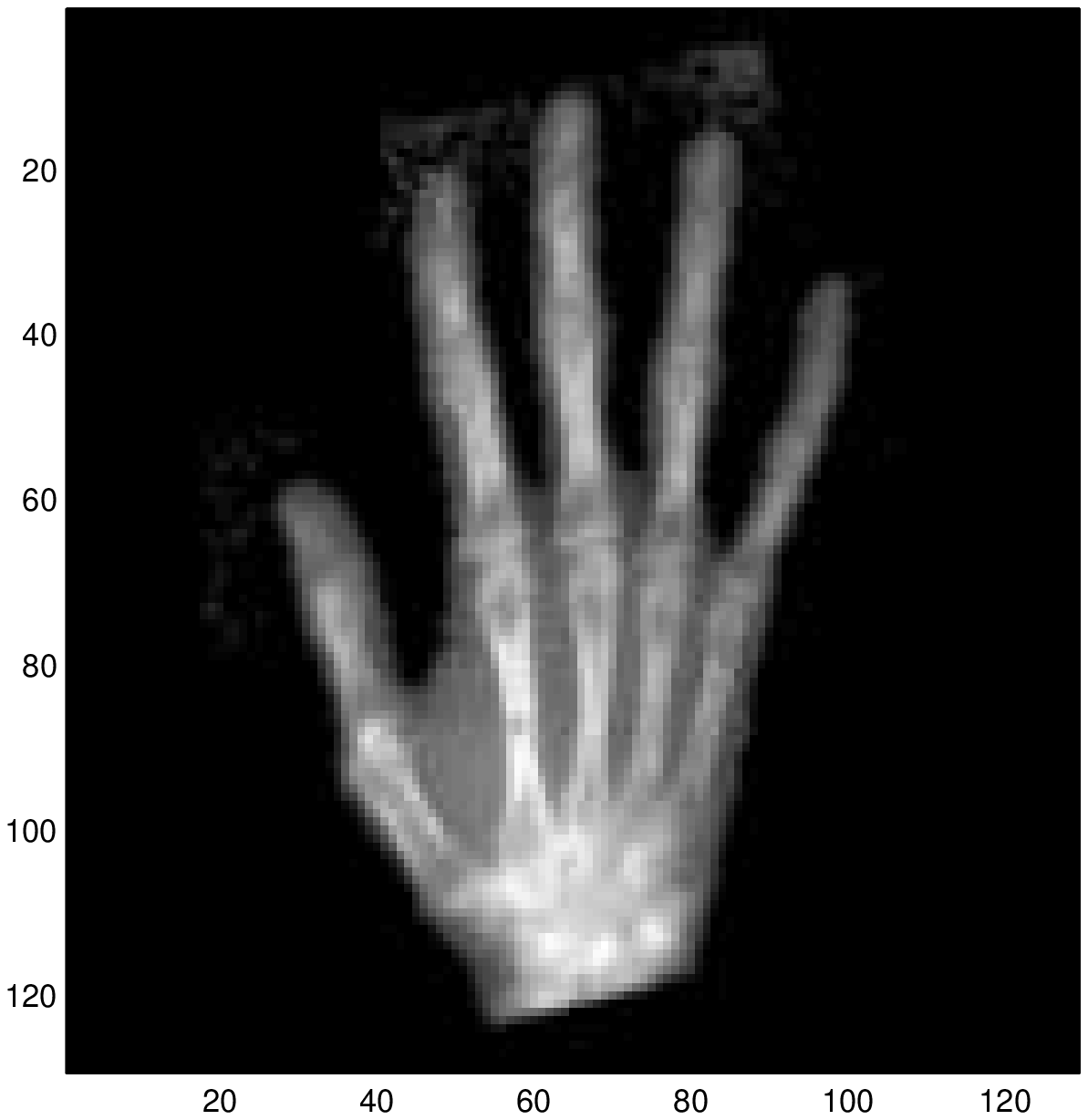}}
    & \subfloat[$T-R$]{\includegraphics[height=4cm,width=4.5cm]{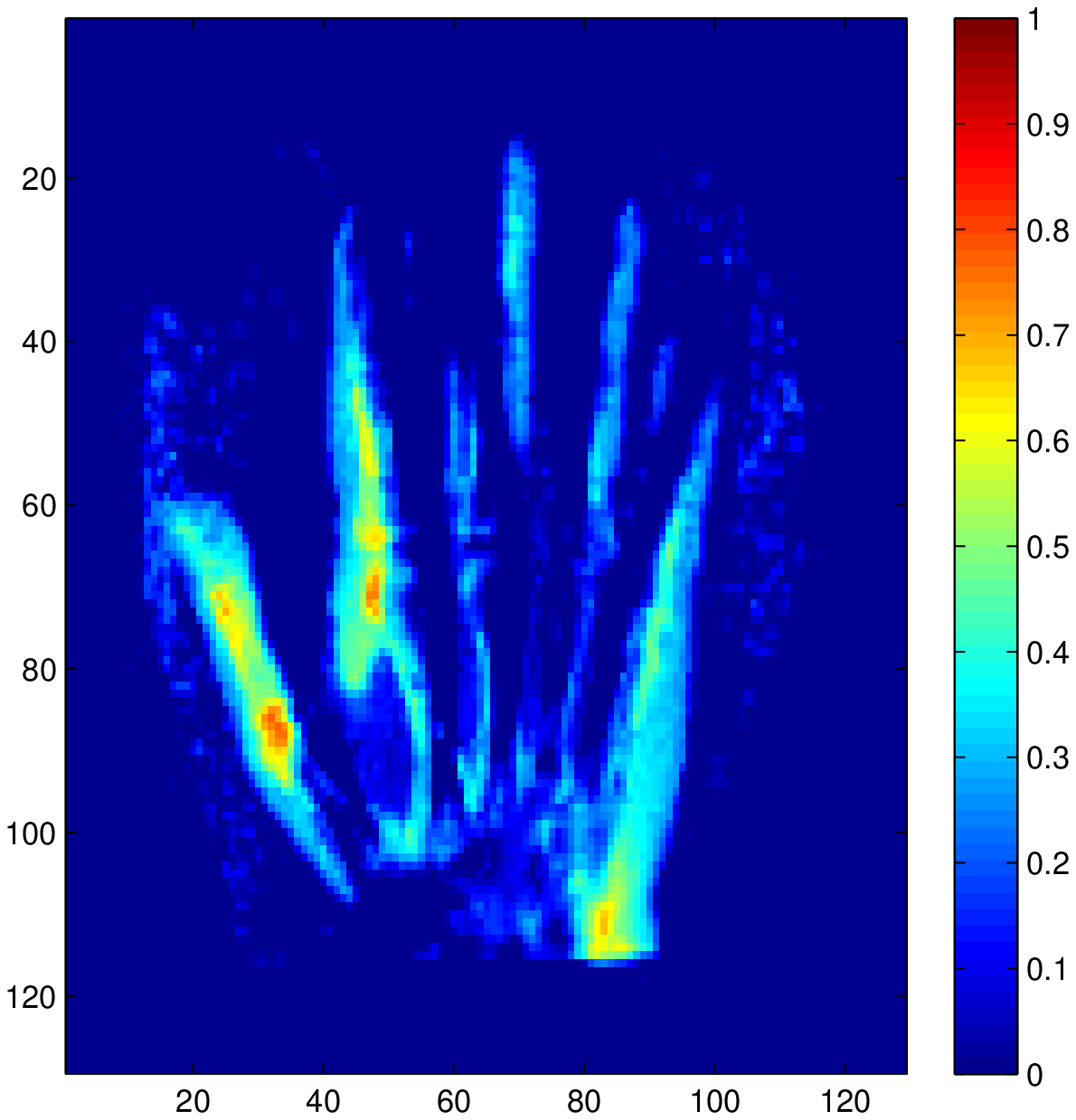}}\\
    \subfloat[$\boldsymbol{x}+\boldsymbol{u}(\boldsymbol{x})$]
    {\includegraphics[height=4cm,width=4.5cm]{G_1.eps}}
   & \subfloat[$T(\boldsymbol{x}+\boldsymbol{u}(\boldsymbol{x}))$]
   {\includegraphics[height=4cm,width=4.5cm]{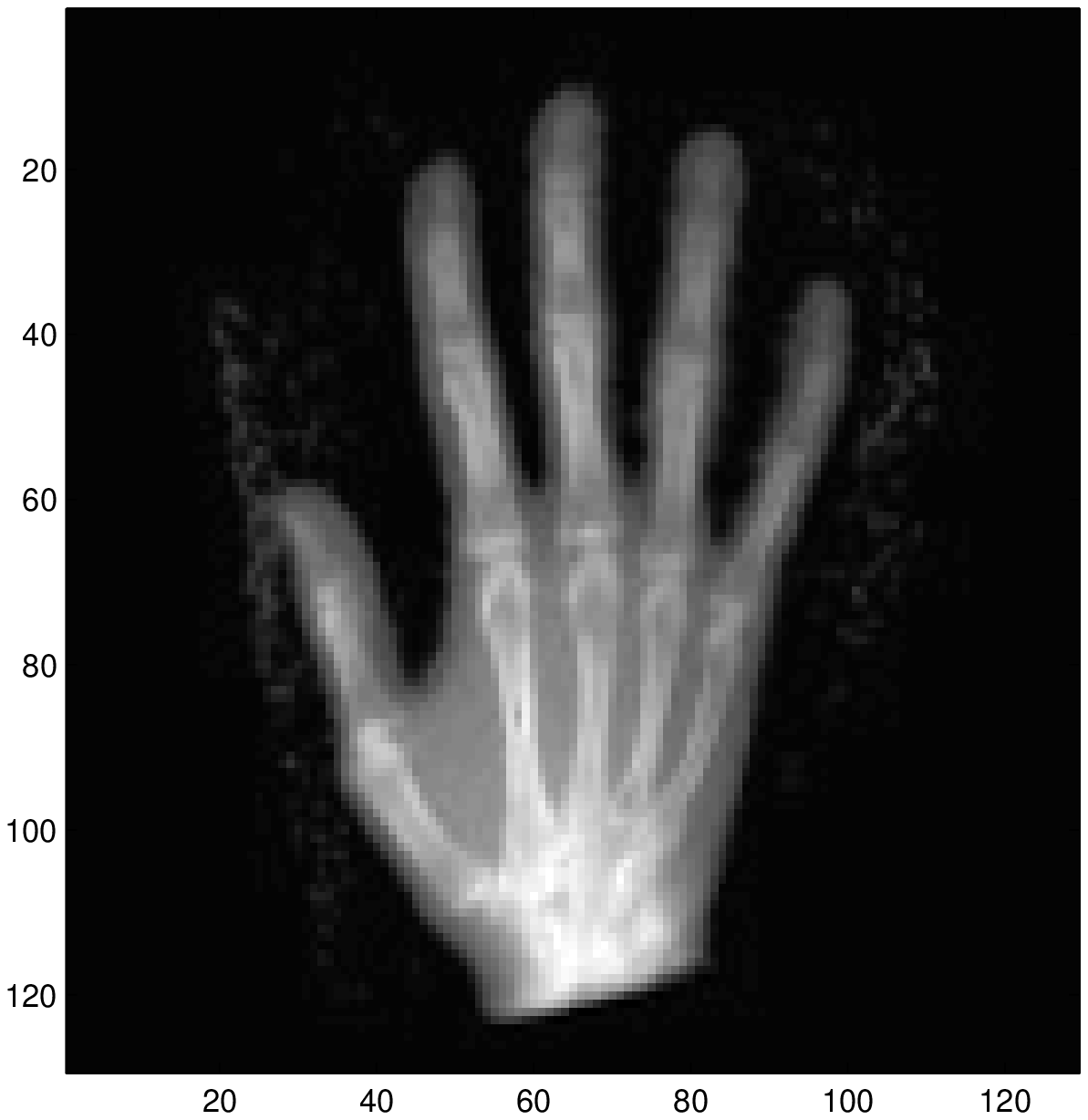}}
    & \subfloat[$T(\boldsymbol{x}+\boldsymbol{u}(\boldsymbol{x}))-R$]
    {\includegraphics[height=4cm,width=4.5cm]{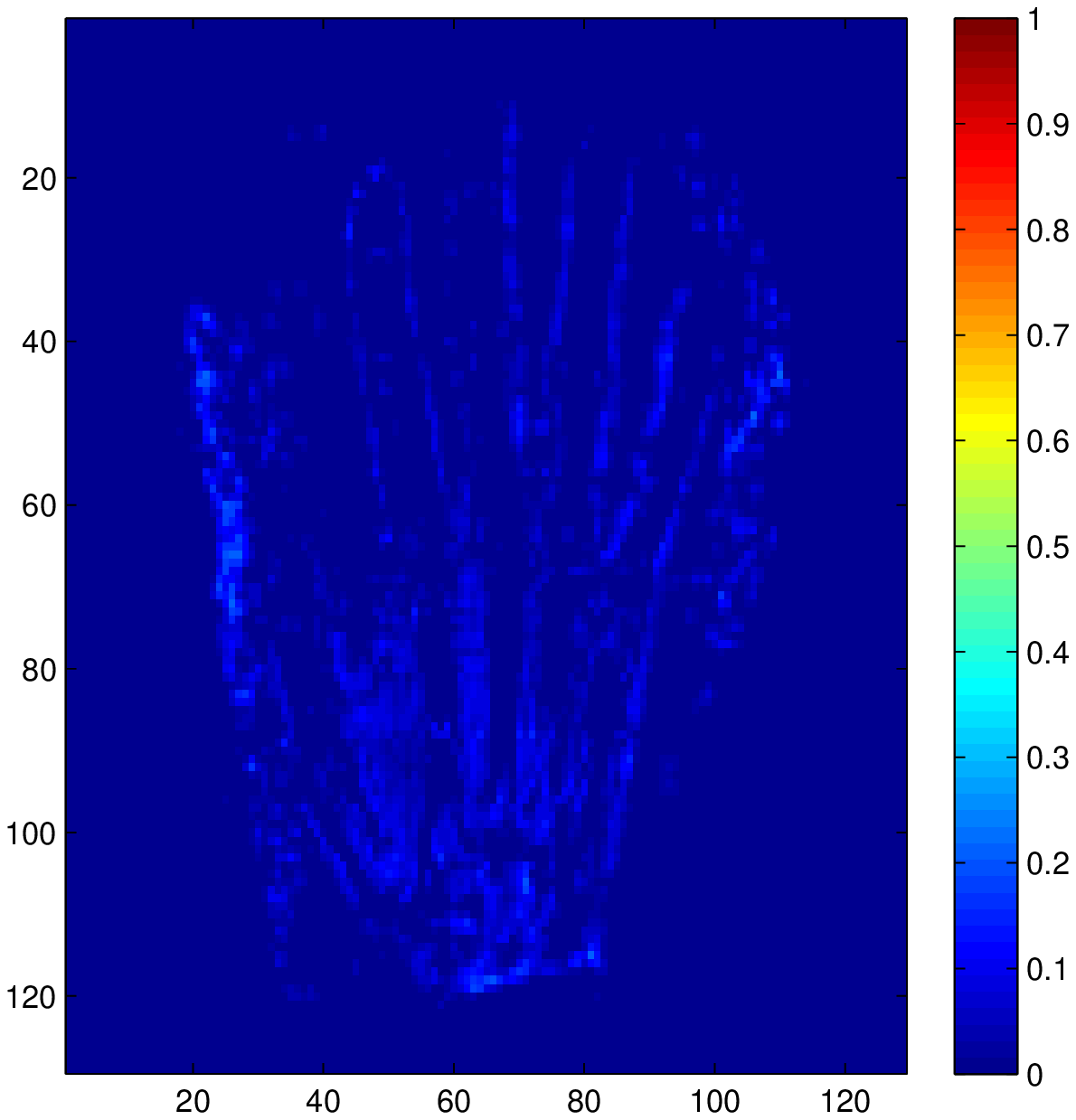}}\\
\end{tabular}
%\caption{Test 1 (X-ray hand). Top row and left to right: template, reference and the difference before registration. Bottom row and left to right: the transformation is applied to a regular grid, the transformed template image and the difference after registration. Gaussian curvature is able to solve smooth problems. }
\caption{Test 1 (X-ray hand). Illustration of the effectiveness of Gaussian curvature with smooth problems. On the top row, from left to right: (a) template, (b) reference and (c) the difference before registration. On the bottom row, from left to right: (d) the transformation applied to a regular grid, (e) the transformed template image and (f) the difference after registration. As can be seen from the result (e) and the small difference after registration (f), Gaussian curvature is able to solve smooth problems. }
\label{f1}
\end{figure}

\begin{figure}[h!]
\centering
\begin{tabular}{cc}%\hspace*{-6mm}
 \subfloat[Model D]{\includegraphics[height=4.5cm,width=5cm]{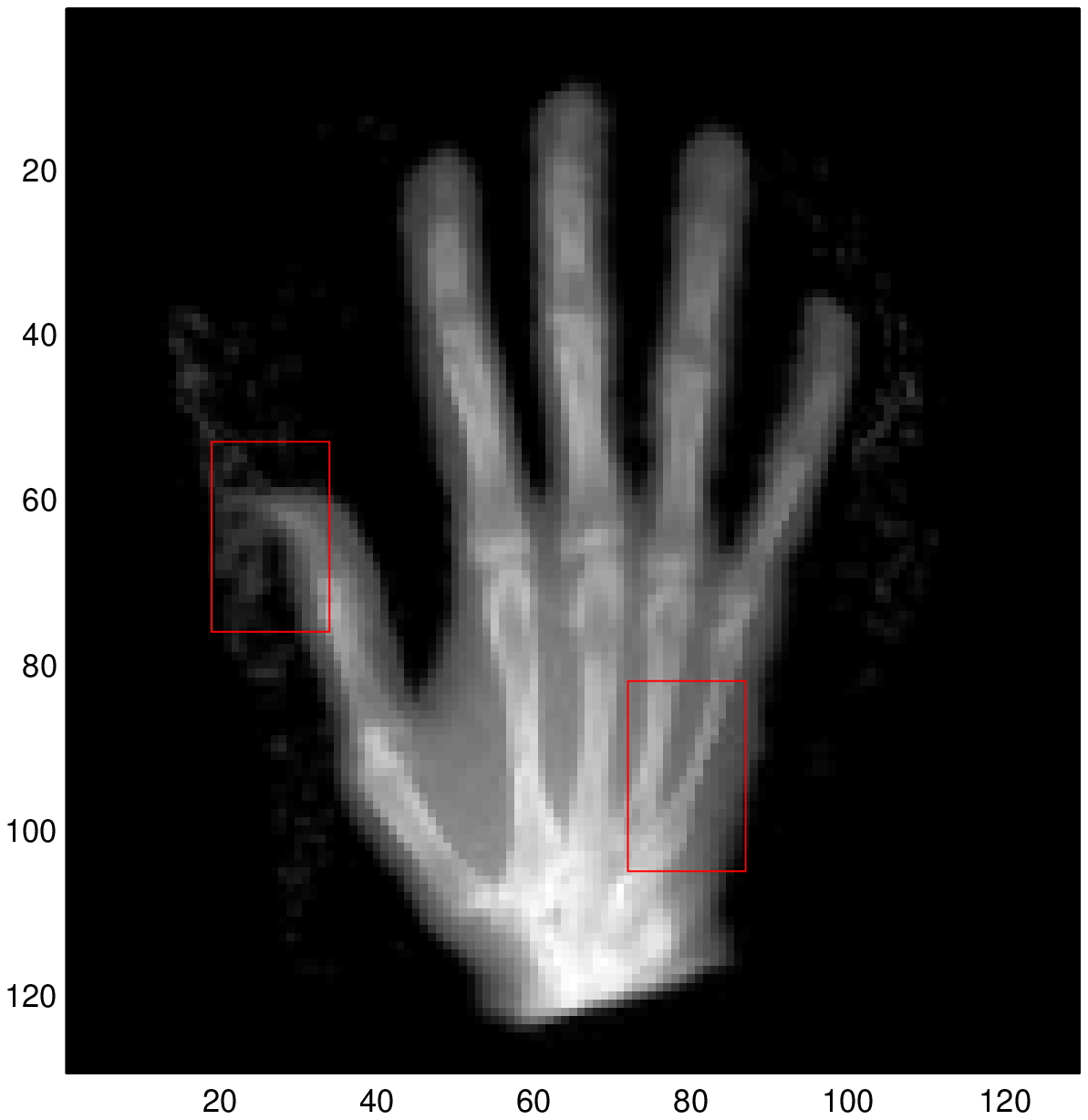}}
   & \subfloat[Model LC]{\includegraphics[height=4.5cm,width=5cm]{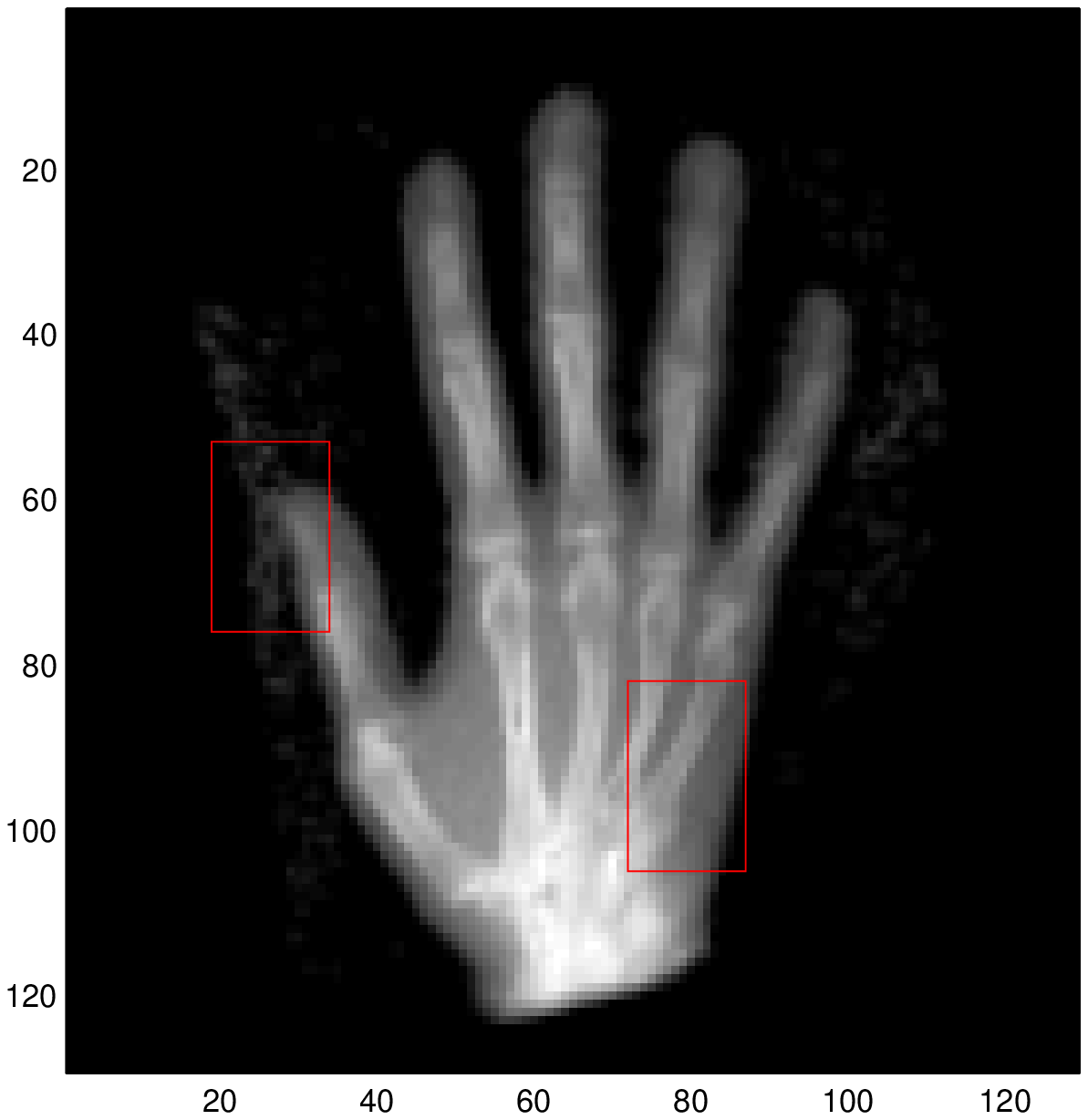}}\\
  \subfloat[Model MC]{\includegraphics[height=4.5cm,width=5cm]{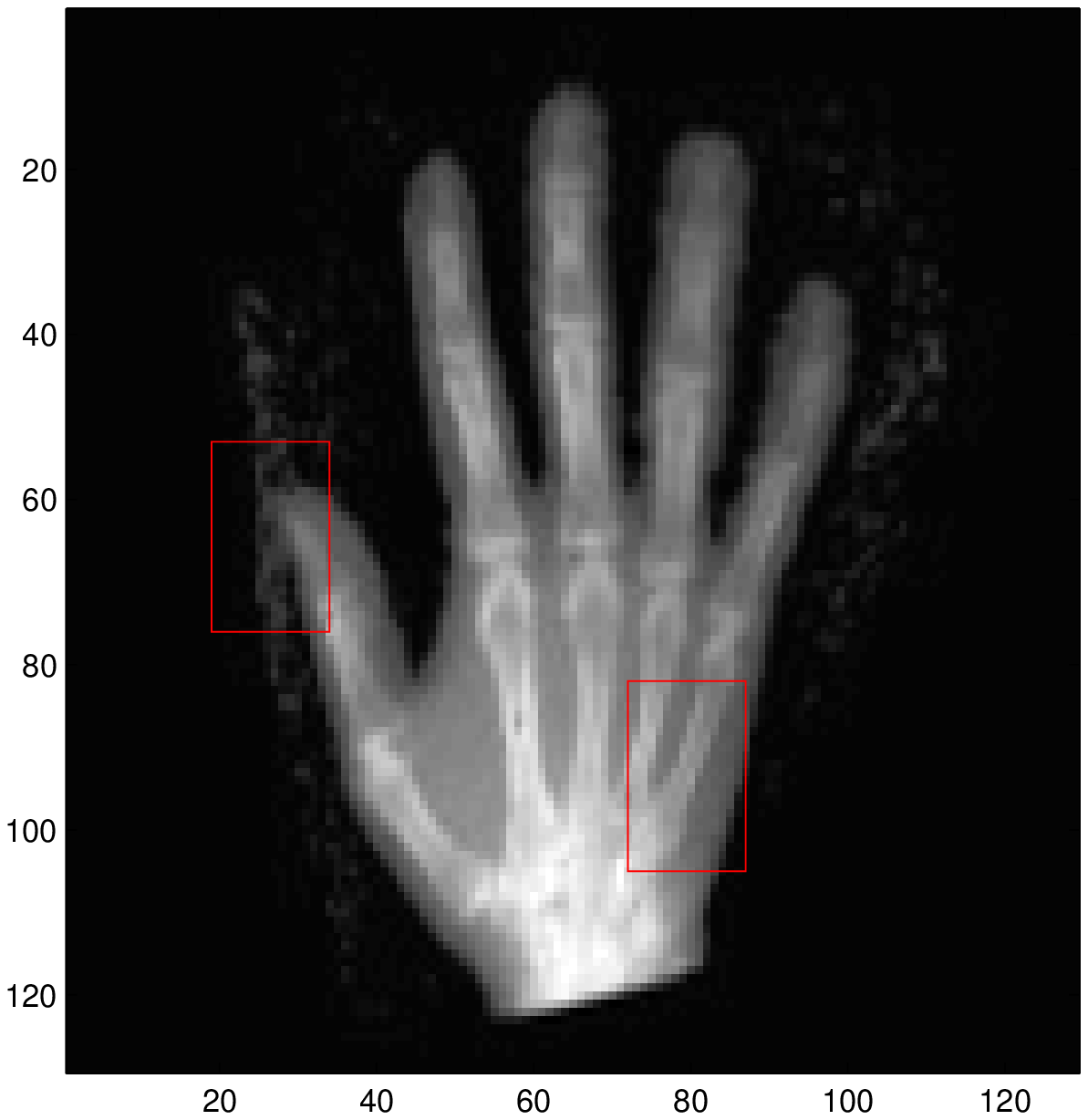}}
    & \subfloat[Gaussian curvature]{\includegraphics[height=4.5cm,width=5cm]{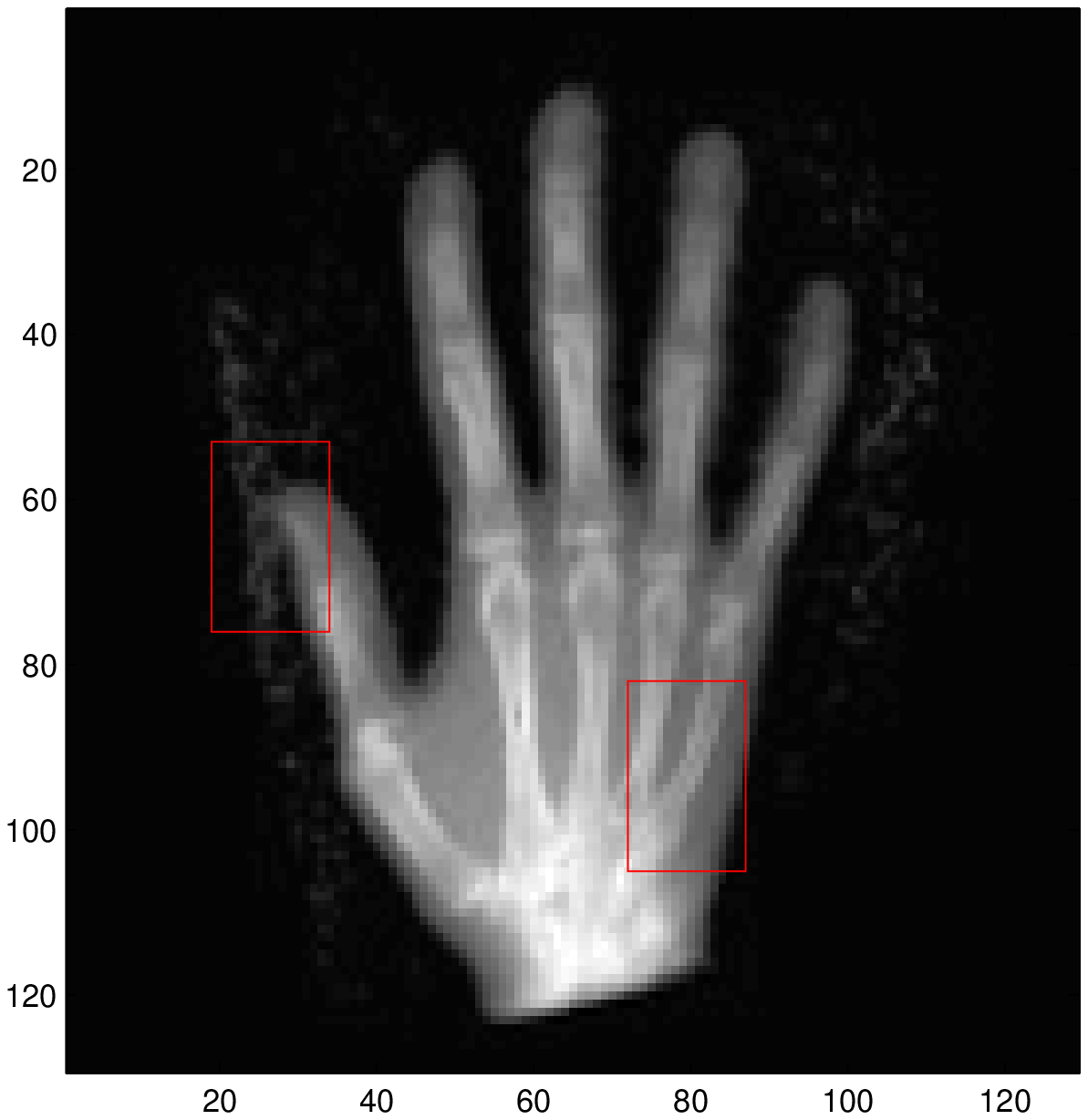}}\\
\end{tabular}
\caption{Test 1 (X-ray hand). Comparison of Gaussian curvature with competing methods. The transformed template image using (a) Model D, (b) Model LC,
(c) Model MC and (d) Gaussian curvature. Note the difference of these three images inside the red boxes.}
\label{f1a}
\end{figure}

Images for Test 1 are taken from \cite{FAIR} where X-ray images of two hands of different individuals need to be aligned. The size of images are $128\times 128$ and the recovered transformation is expected to be smooth.  The scaled version of the transformation  and the transformed template image is given in Figures \ref{f1} (d) and (e) respectively. The transformation is smooth and the model is able to solve such a problem.
%%%%%%%Result for other method
For comparison, the transformed template images for the diffeormorphic demon method, linear, mean and Gaussian curvatures are shown in Figures \ref{f1a} (a), (b), (c) and (d) respectively. We can observe that there are some differences of these images inside the red boxes where only Gaussian curvature delivering the best result of the features inside the boxes.

%%%table 1 for summarize the results
\begin{table}[h!]
\centering
\def\tabularxcolumn#1{m{#1}}
\begin{tabular}{|c|c|c|c|c|c|c|} \hline
Measure &  \multicolumn{2}{c|}{Model D} & \multicolumn{2}{c|}{Model LC} & Model MC & GC\\  \hline
{}      & ML   & SL                 & ML   & SL             &SL & SL \\  \hline
$\gamma$&  1.6 & 1.6               & 0.1   & 0.5            & 0.0001 & 0.0001\\ \hline
Time (s)&  15.19& 186.48               & 84.33   & 12.98           & 275.3   & 953.15 \\ \hline
$\varepsilon$&  0.1389 & 0.1229               & 0.0720  & 0.3780            & 0.0964  &0.0582  \\ \hline
${\cal F}$&  0.0600 & 0.1082              & 0.3894   & 0.1973            & 0.6390 & 0.3264\\ \hline
\end{tabular}
\caption{Quantitative measurements for all models for Test 1. ML and SL stand for multi and single level respectively. $\gamma$ is chosen as small as possible such that ${\cal F}>0$ for all methods. ${\cal F}>0$ indicates the deformation consists of no folding and cracking of the deformed grid. We can see that the smallest value of $\varepsilon$ is given by Gaussian curvature (GC). }\label{t1}
\end{table}

We summarised the results for Test 1 in Table \ref{t1} where ML and SL stand for multi and single level respectively. For all models, $\gamma$ is chosen as small as possible such that ${\cal F}>0$. We can see that the fastest model is the diffeormorphic demon, followed by linear and mean curvature. The current implementation for Gaussian curvature is on single level and the model uses augmented Lagrangian method which has four dual variables and four lagrange multipliers terms. Thus, it requires more computational time than the other models.

\subsection{Test 2: A Pair of Brain MR Images}
\begin{figure}[h!]
\centering
\def\tabularxcolumn#1{m{#1}}
\begin{tabular}{ccc}%\hspace*{-6mm}
 \subfloat[$T$]{\includegraphics[height=4cm,width=4.5cm]{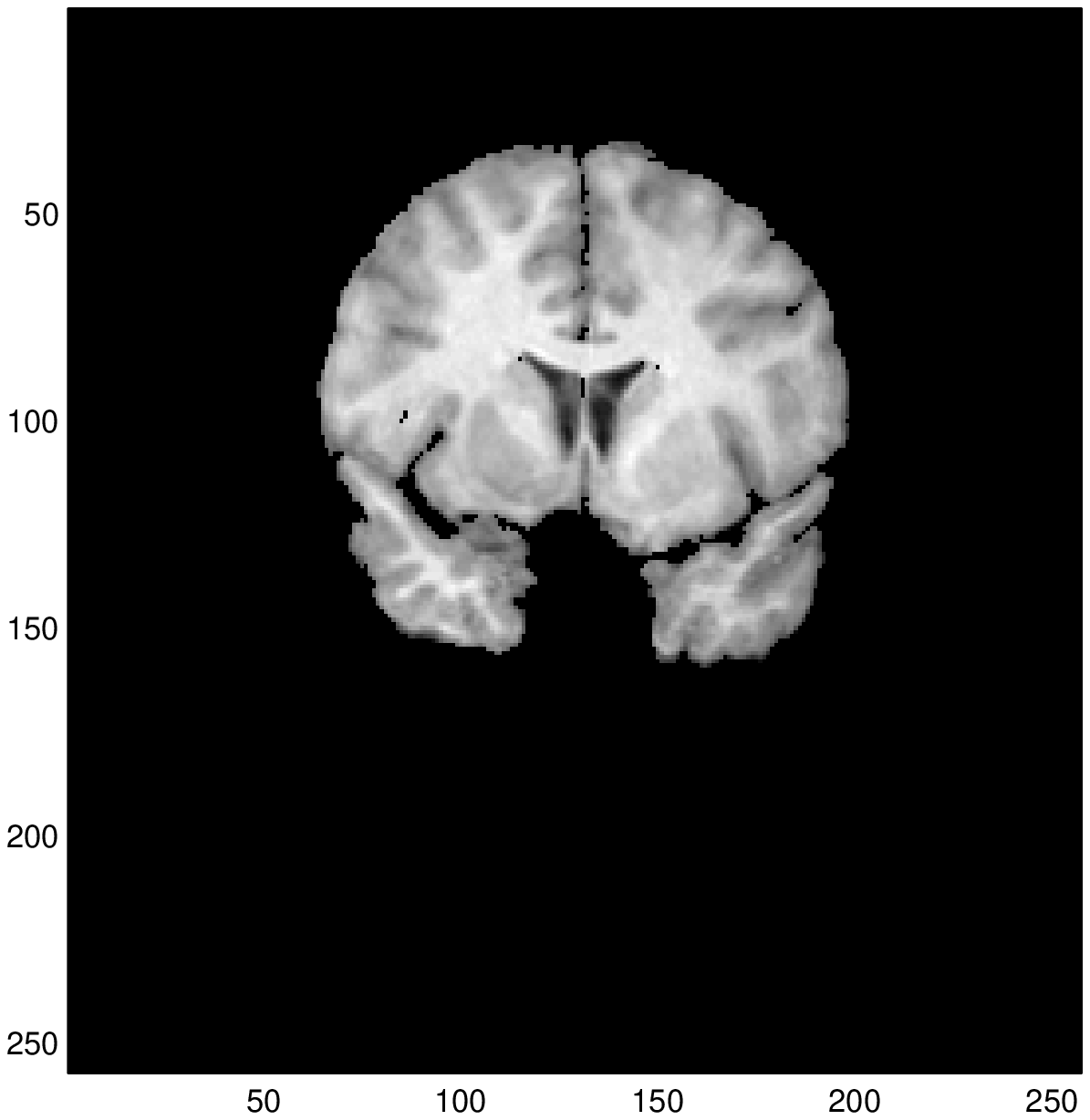}}
   & \subfloat[$R$]{\includegraphics[height=4cm,width=4.5cm]{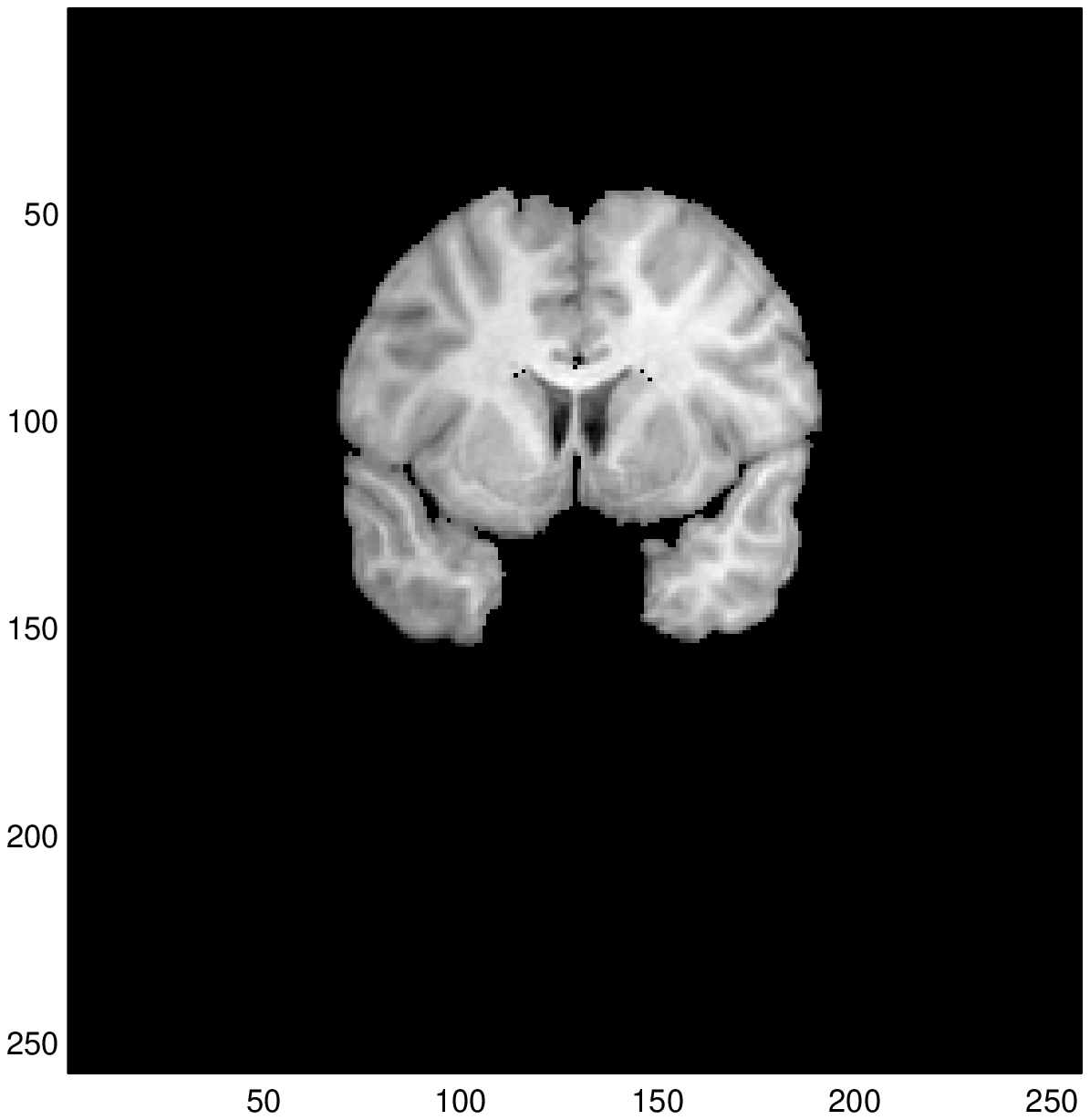}}
    & \subfloat[$T-R$]{\includegraphics[height=4cm,width=4.5cm]{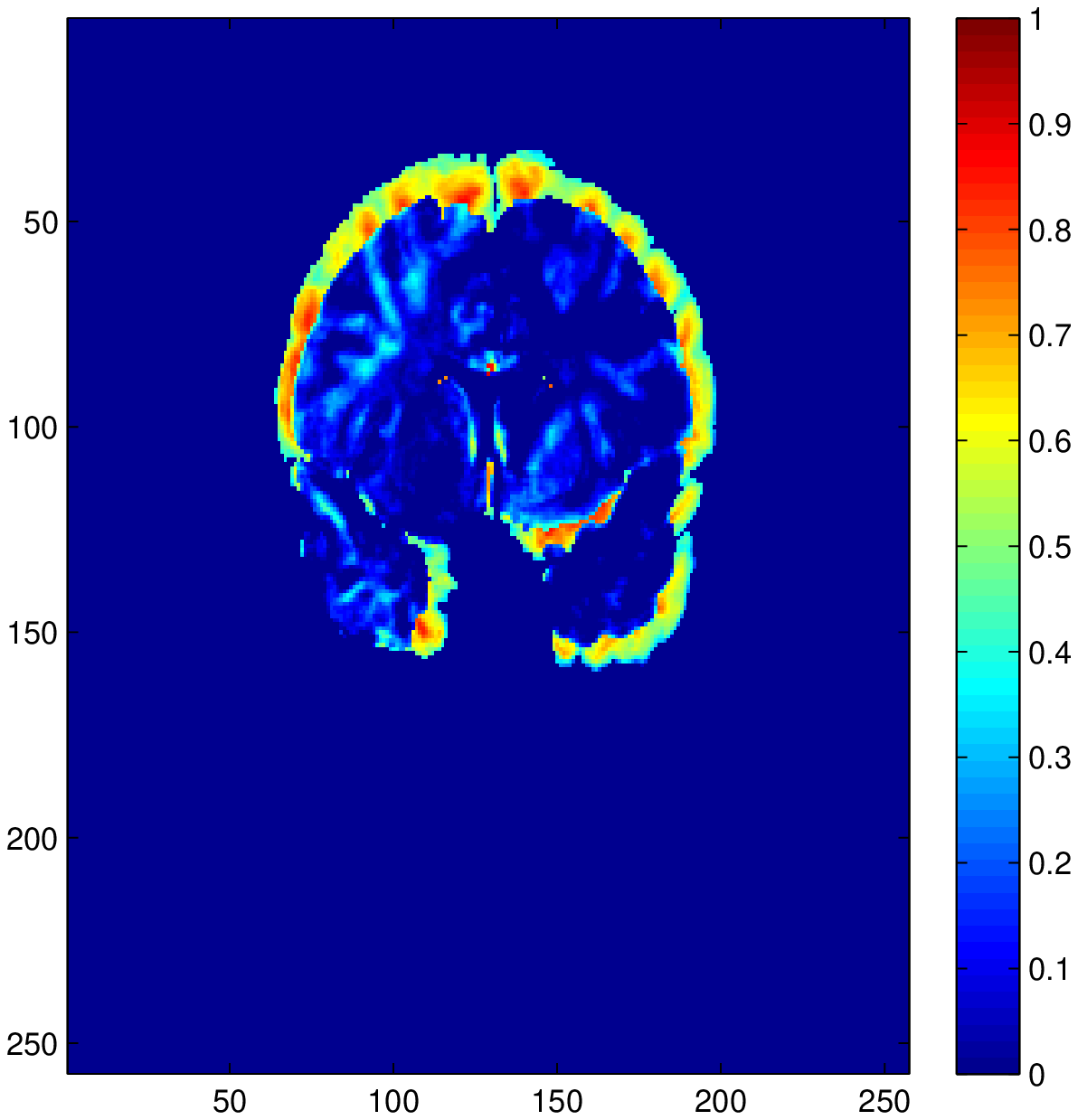}}\\
    \subfloat[$\boldsymbol{x}+\boldsymbol{u}(\boldsymbol{x})$]
    {\includegraphics[height=4cm,width=4.5cm]{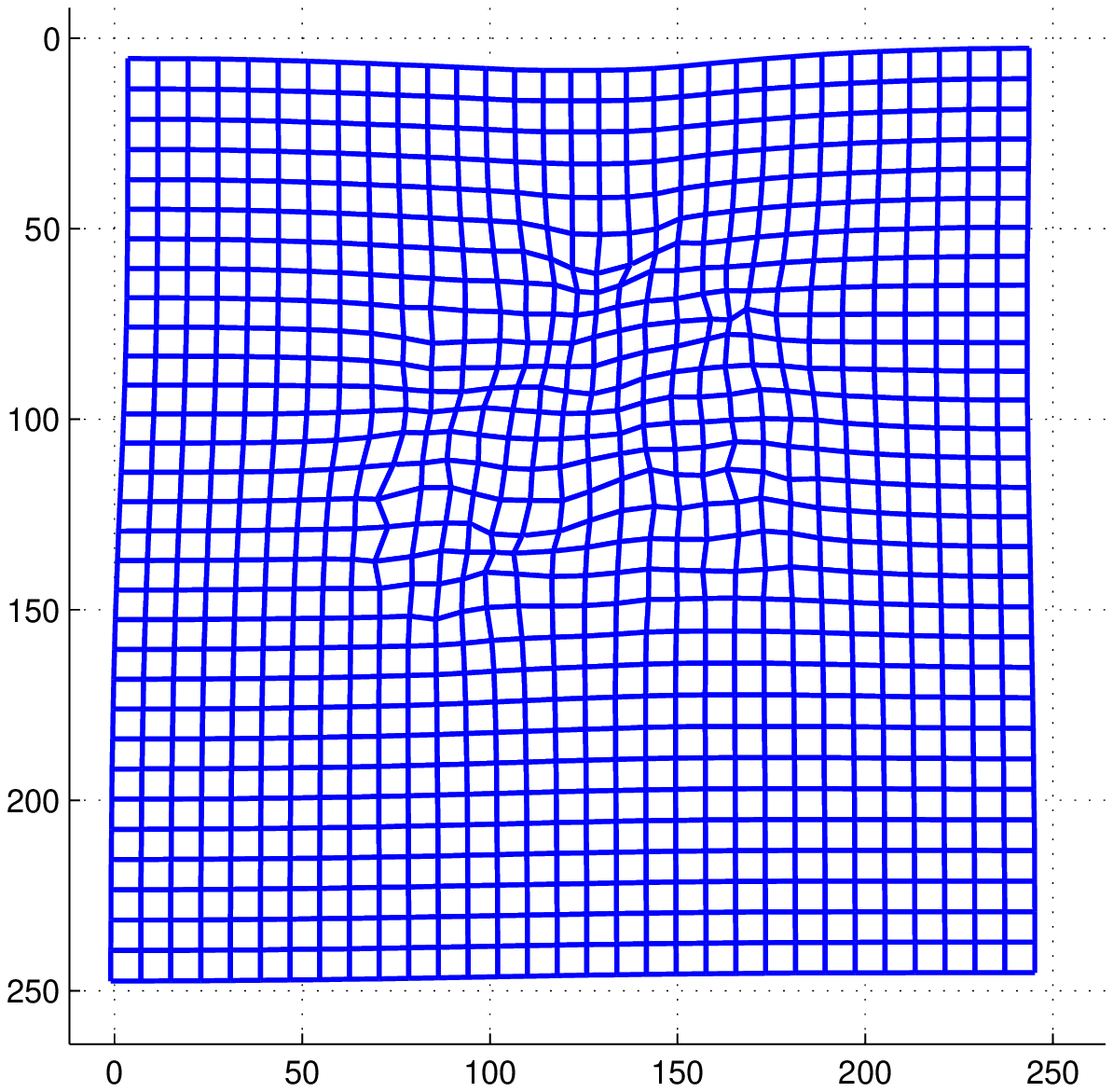}}
   & \subfloat[$T(\boldsymbol{x}+\boldsymbol{u}(\boldsymbol{x}))$]
   {\includegraphics[height=4cm,width=4.5cm]{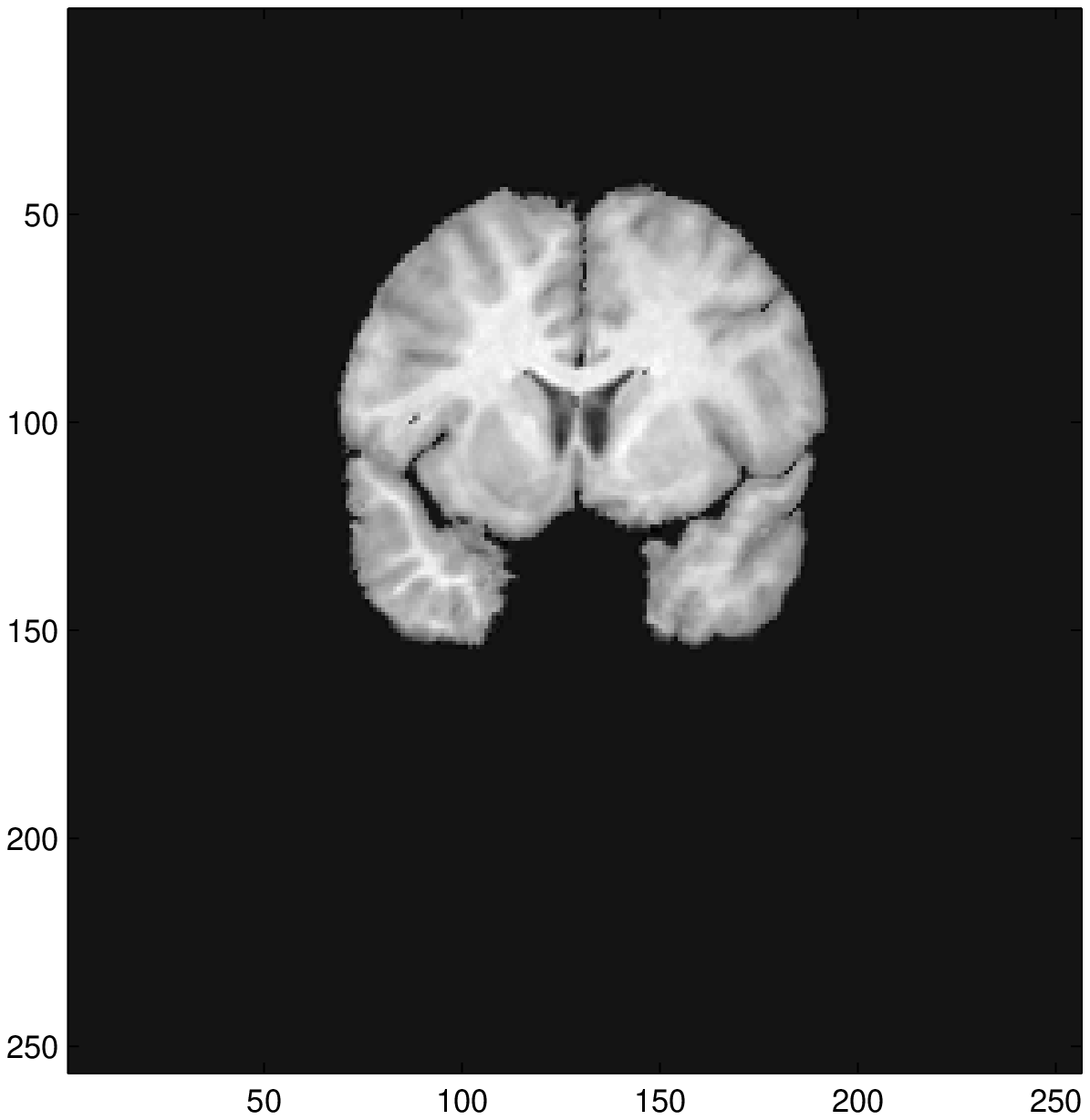}}
    & \subfloat[$T(\boldsymbol{x}+\boldsymbol{u}(\boldsymbol{x}))-R$]
    {\includegraphics[height=4cm,width=4.5cm]{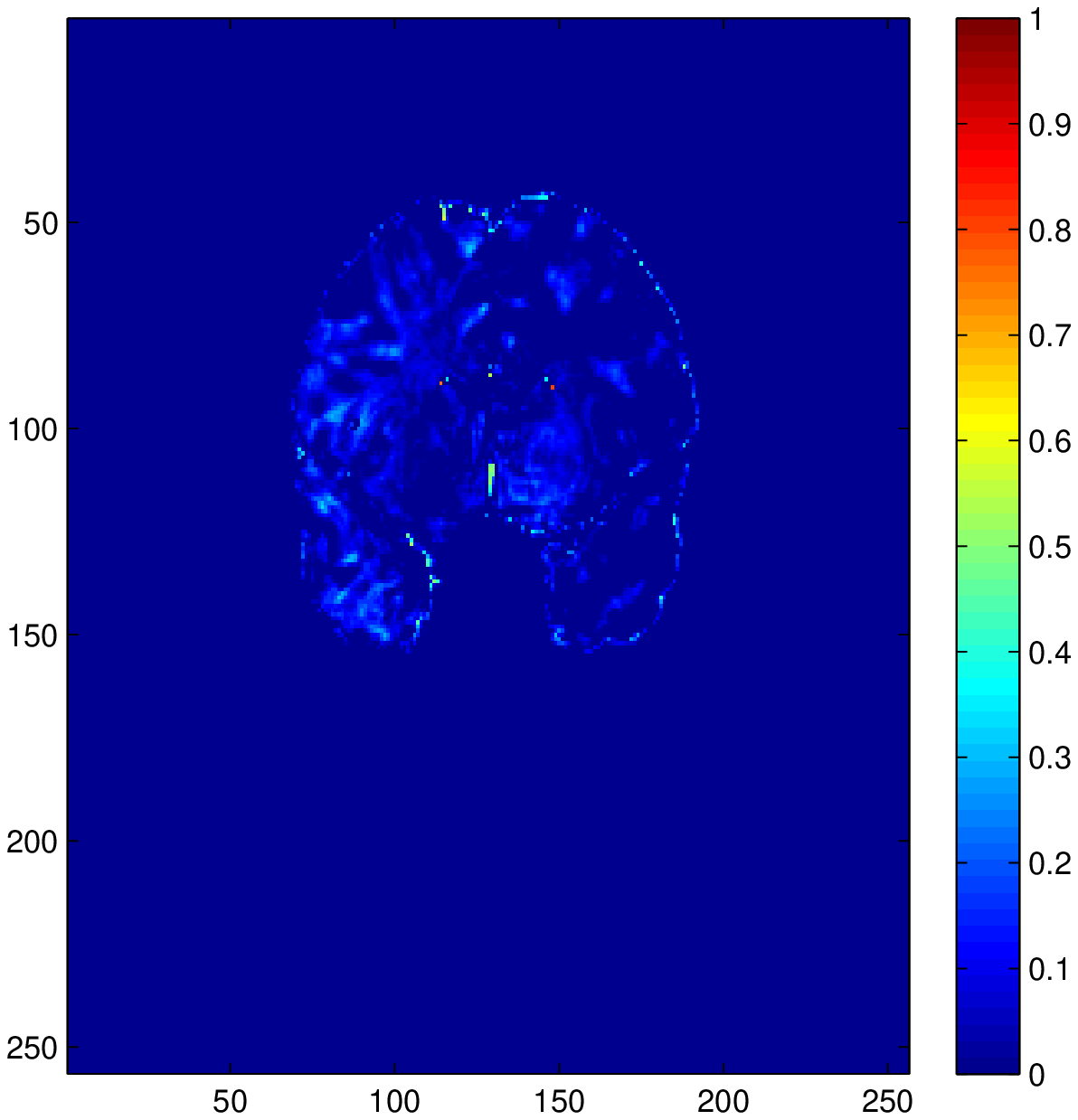}}\\
\end{tabular}
%\caption{Test 2: A pair of Brain MR images. Top row and left to right: template, reference and the difference before registration. Bottom row and left to right: the transformation is applied to a regular grid, the transformed template image and the difference after registration. Gaussian curvature can be applied to real medical images. }
\caption{Test 2: A pair of Brain MR images. Illustration of the effectiveness of Gaussian curvature with real medical images. On the top row, from left to right: (a) template, (b) reference and (c) the difference before registration. On the bottom row, from left to right: (d) the transformation applied to a regular grid, (e) the transformed template image and (f) the difference after registration. As can be seen from the result (e) and the small difference after registration (f), Gaussian curvature can be applied to real medical images and is able to obtain good results. }
\label{f3}
\end{figure}

\begin{figure}[h!]
\centering
\def\tabularxcolumn#1{m{#1}}
\begin{tabular}{cc}%\hspace*{-6mm}
 \subfloat[Model D]{\includegraphics[height=4cm,width=4.5cm]{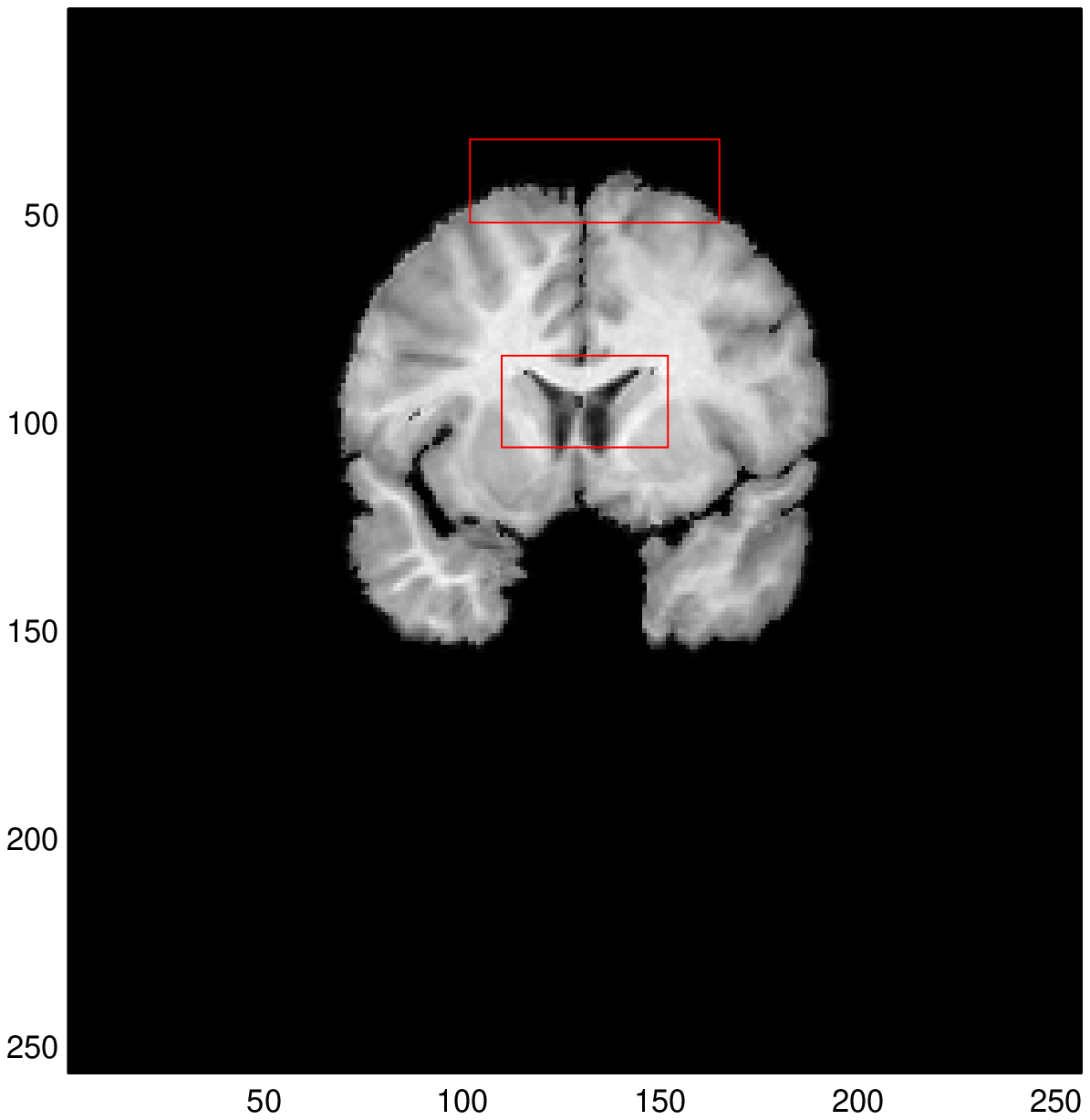}}
   & \subfloat[Model LC]{\includegraphics[height=4cm,width=4.5cm]{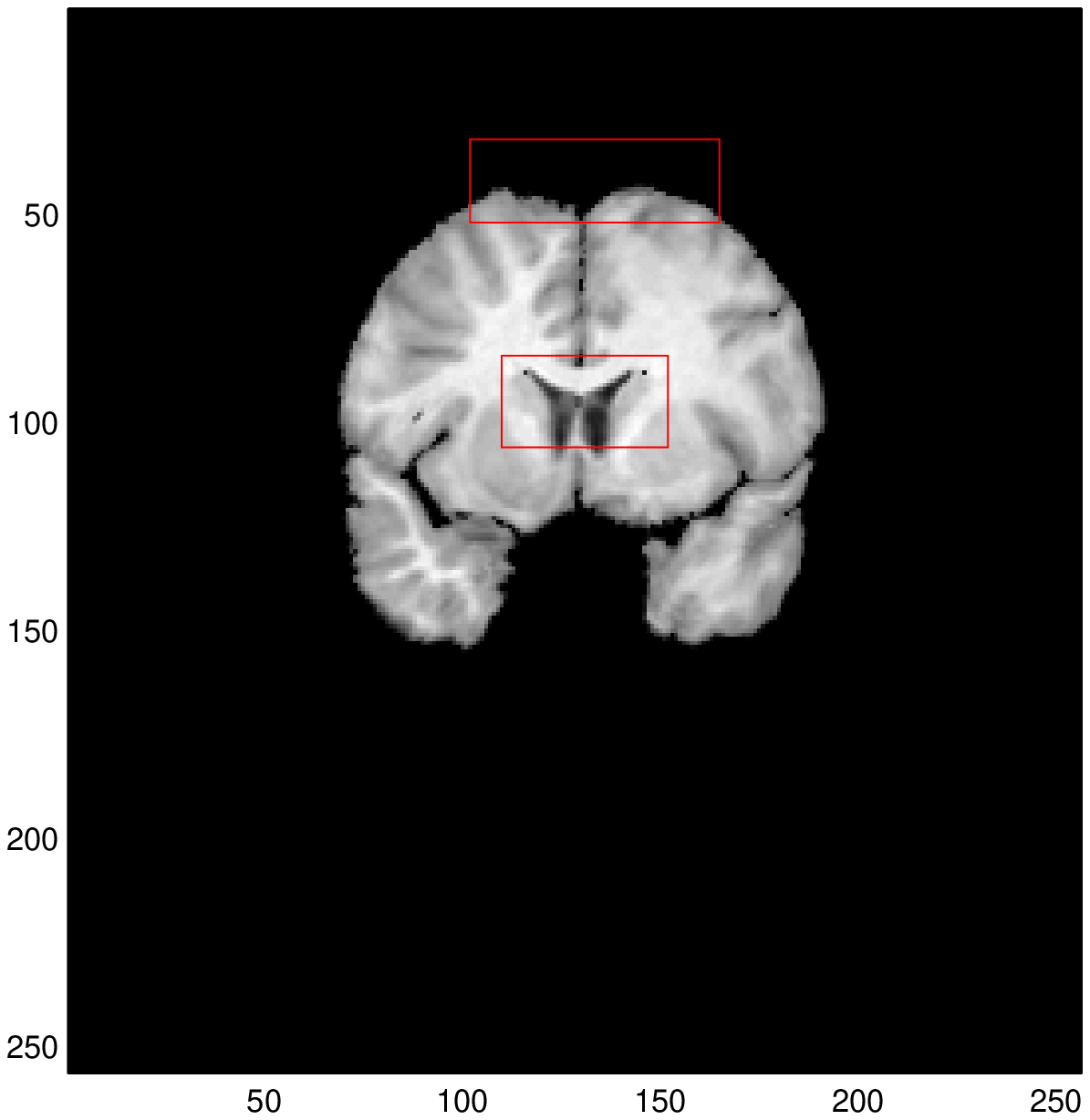}}\\
 \subfloat[Model MC]{\includegraphics[height=4cm,width=4.5cm]{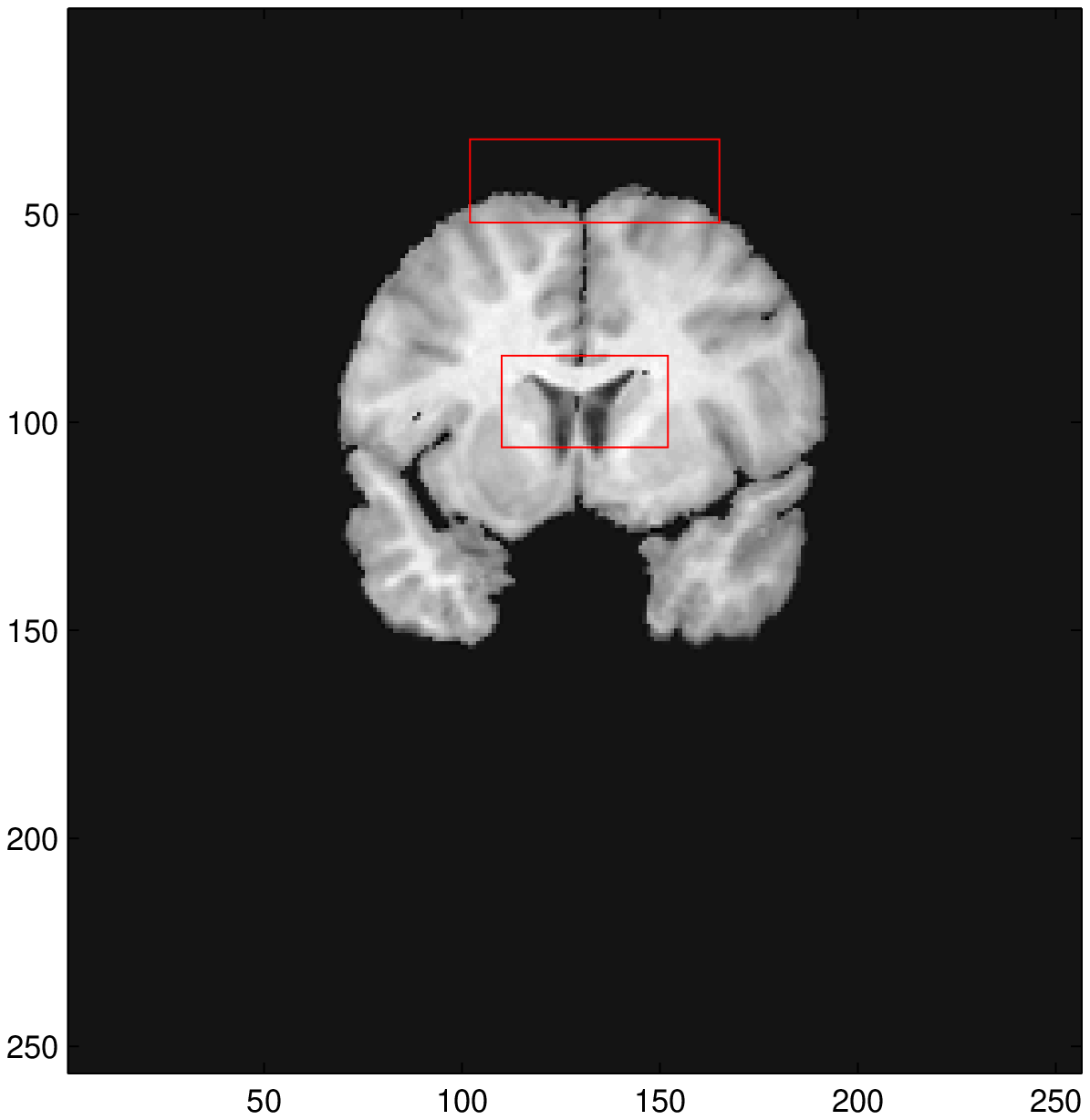}}
    & \subfloat[Gaussian curvature]{\includegraphics[height=4cm,width=4.5cm]{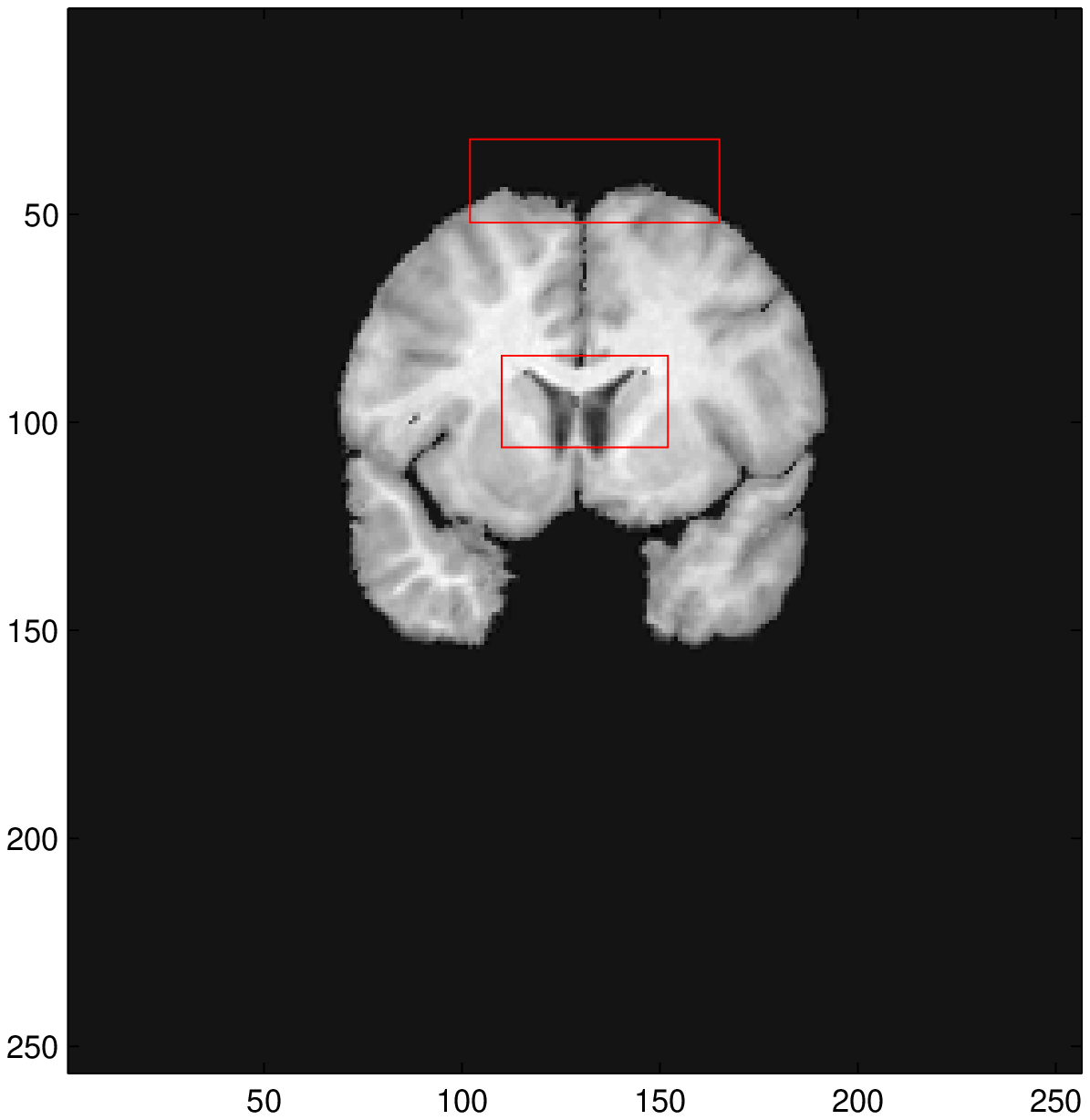}}\\
\end{tabular}
%\caption{Test 2 Left to right. The transformed template image using the diffeomorphic demon, linear and Gaussian curvature. Notice the differences of these three images inside the red boxes. Gaussian curvature gives the best result for this part of the image.}
\caption{Test 2: A pair of Brain MR images. Comparison of Gaussian curvature with competing methods. The transformed template image using (a) Model D,
(b) Model LC, (c) Model MC, and (d) Gaussian curvature. Notice the differences of these three images inside the red boxes. Considerably more accurate results are obtained, particularly within these significant regions, by employment of the Gaussian curvature model.}
\label{f3a}
\end{figure}

We take  as Test 2 a pair of medical images of size $256 \times 256$   from the Internet Brain Segmentation Repository (IBSR)
$\mbox{\tt http://www.cma.mgh.harvard.edu/ibsr} $
where 20 normal MR brain images and their manual segmentations are provided. We choose the particular pair of   individuals with different sizes of ventricle to illustrate a large deformation problem. Figure \ref{f3} shows the test images and the registration results using Gaussian curvature model. We can see that the model is able to solve real medical problems involving large deformations, which is particularly important for atlas construction in medical applications.
%%%%%%%Result for other method
Figure \ref{f3a} shows the transformed template images for all four methods. We can see that Gaussian curvature gives the best result inside the red boxes in comparison with the diffeomorphic demon, the linear and mean curvature models as depicted in Figure \ref{f3a} (d).

%\begin{figure}[h!]
%\centering
%\def\tabularxcolumn#1{m{#1}}
%\begin{tabular}{|c|c|c|c|c|}
%  \hline
%  Test & Diff. Demon & Linear Curv. & Mean Curv. &Gaussian Curv.\\ \hline
%  1 & 0.1229  & 0.0836 & 0.1060 &0.0663\\ \hline
%2  & 0.2245  & 0.1914  & 0.2098      &0.1311\\
%  \hline
%\end{tabular}
%\caption{Values of the reduction-of-similarity measure $\varepsilon$ for Tests 1 and 2 for all of the methods. We can see that the lowest values of $\varepsilon$ are given by Gaussian curvature for both tests indicating higher measure of alignment between the resulting transformed template image and the original reference image. }\label{t1}
%\end{figure}

\begin{table}[h!]
\centering
\def\tabularxcolumn#1{m{#1}}
\begin{tabular}{|c|c|c|c|c|c|c|} \hline
Measure &  \multicolumn{2}{c|}{Model D} & \multicolumn{2}{c|}{Model LC} & Model MC & GC \\  \hline
{}      & ML   & SL                 & ML   & SL             &SL & SL \\  \hline
$\gamma$&  1.2 & 1.4               & 0.16   & 2.0           & 0.0001 & 0.0001\\ \hline
Time (s)&  23.89& 209.00              & 275.04   & 35.70           & 830.22 &  1053.7\\ \hline
$\varepsilon$&  0.2004& 0.7580             & 0.1128 & 0.4283            & 0.1998&  0.1062 \\ \hline
${\cal F}$&  0.0277& 0.0387             & 0.3157  & 0.0148           & 0.8240 & 0.0138\\ \hline
\end{tabular}
\caption{Quantitative measurements for all models for Test 2. ML and SL stand for multi and single level respectively. $\gamma$ is chosen to be as small as possible such that ${\cal F}>0$ for all models. ${\cal F}>0$ indicates the deformation consists of no folding and cracking of the deformed grid. We can see that the smallest value of $\varepsilon$ is given by Gaussian curvature (GC). }\label{t2}
\end{table}

The values of the quantitative measurements for Test 2 are recorded in Table \ref{t2} where the lowest values of $\varepsilon$ are given by the Gaussian curvature model indicating higher similarity between the transformed template result and the reference image. However, our proposed model required more time than the other models since the model consists more variables than the others.

\begin{figure}[h!]
\centering
\begin{tabular}{cc}%\hspace*{-6mm}
 \subfloat[$T$]{\includegraphics[height=4.5cm,width=5cm]{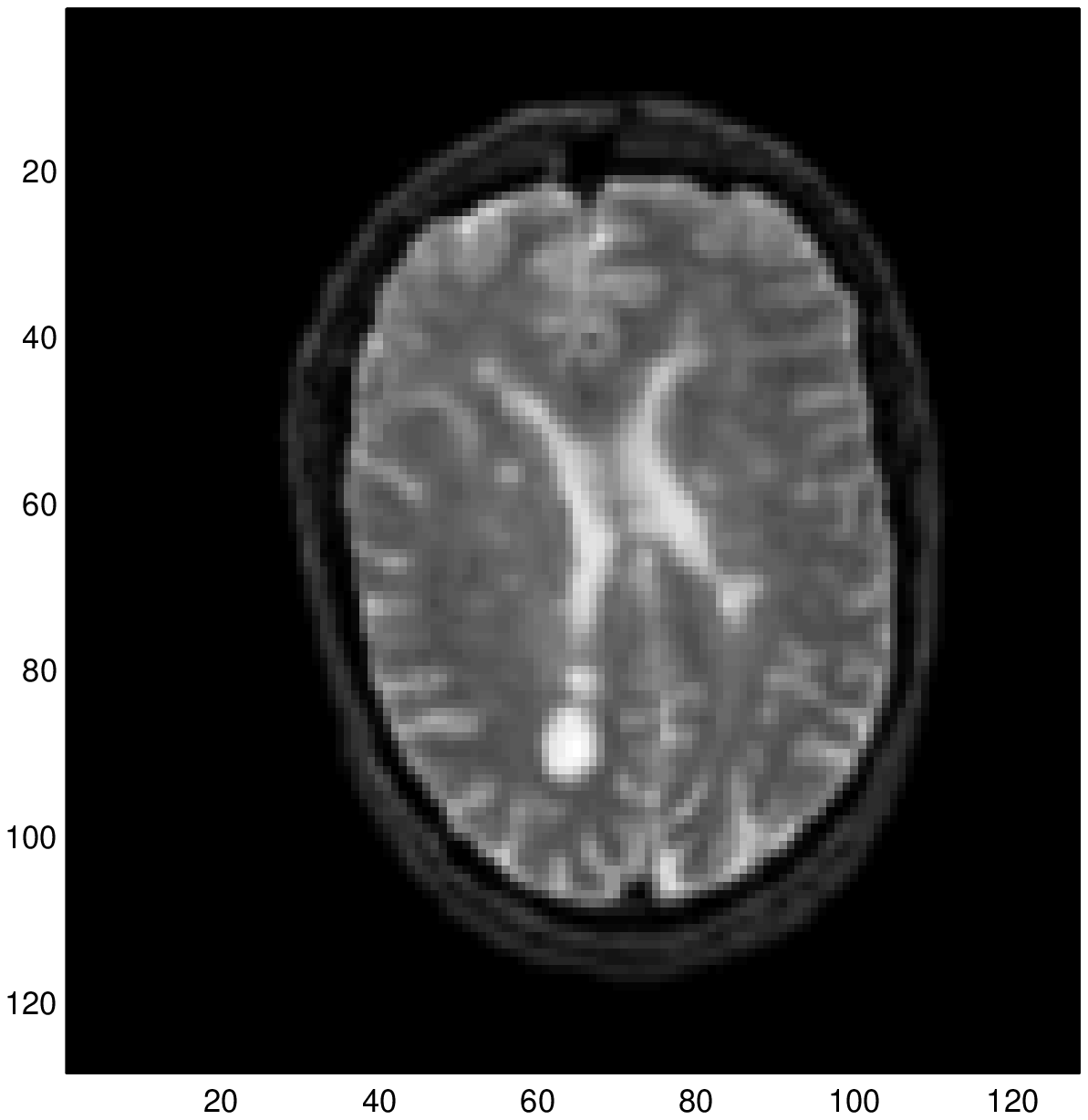}}
   & \subfloat[$R$]{\includegraphics[height=4.5cm,width=5cm]{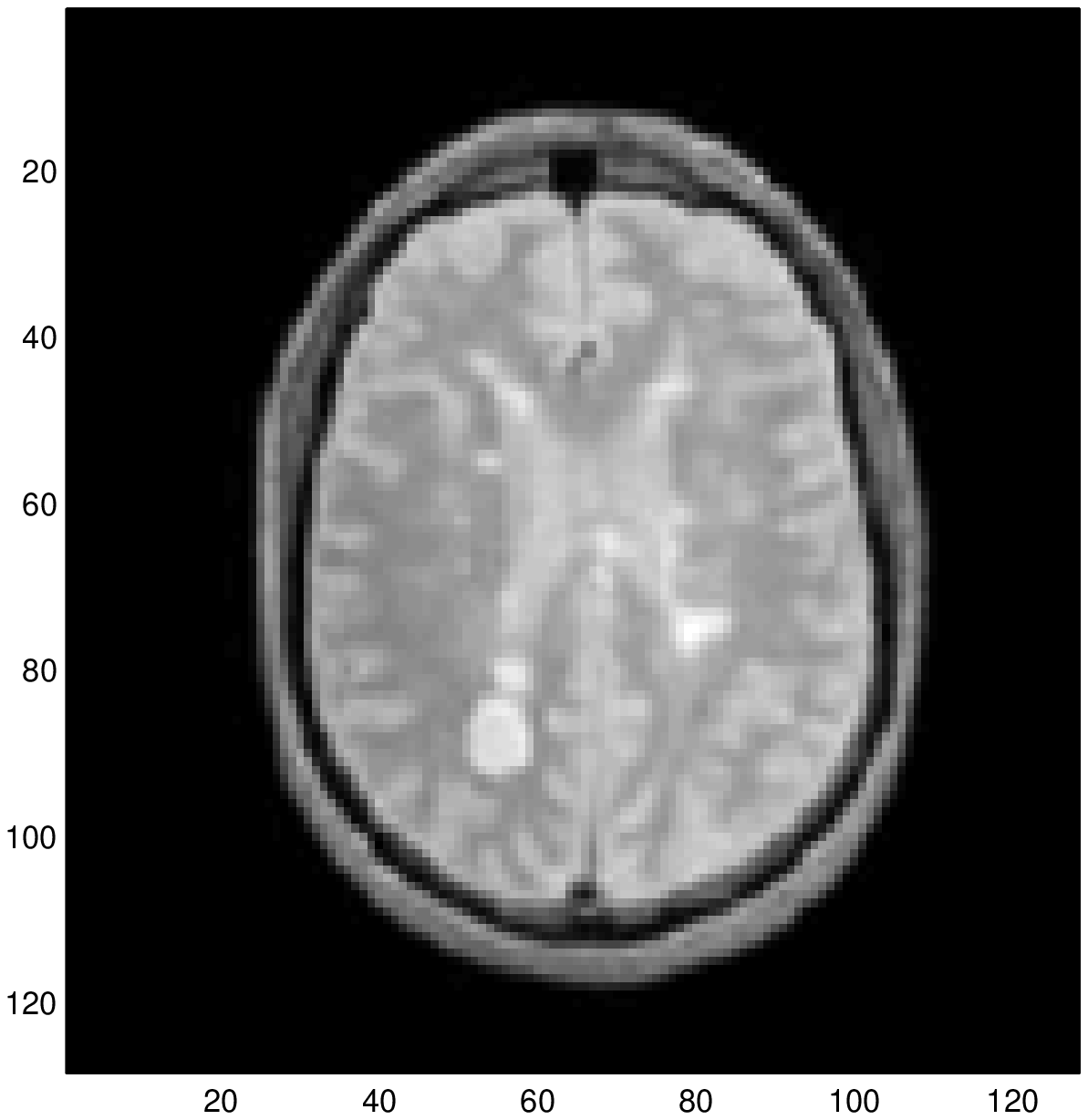}}\\
     \subfloat[$T(\boldsymbol{x}+\boldsymbol{u}(\boldsymbol{x}))$]
     {\includegraphics[height=4.5cm,width=5cm]{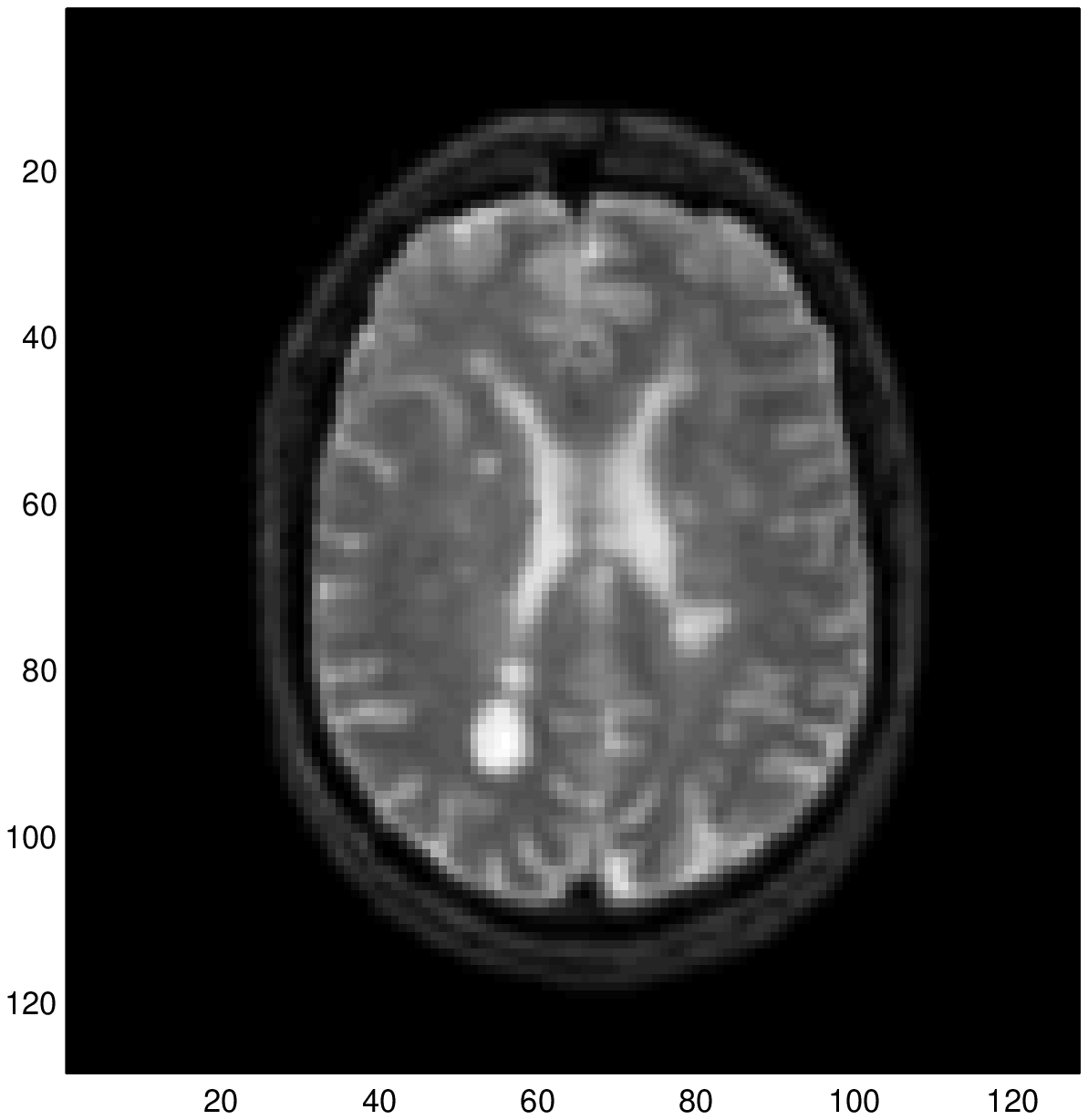}}
    &  \subfloat[$\boldsymbol{x}+\boldsymbol{u}(\boldsymbol{x})$]
    {\includegraphics[height=4.5cm,width=5cm]{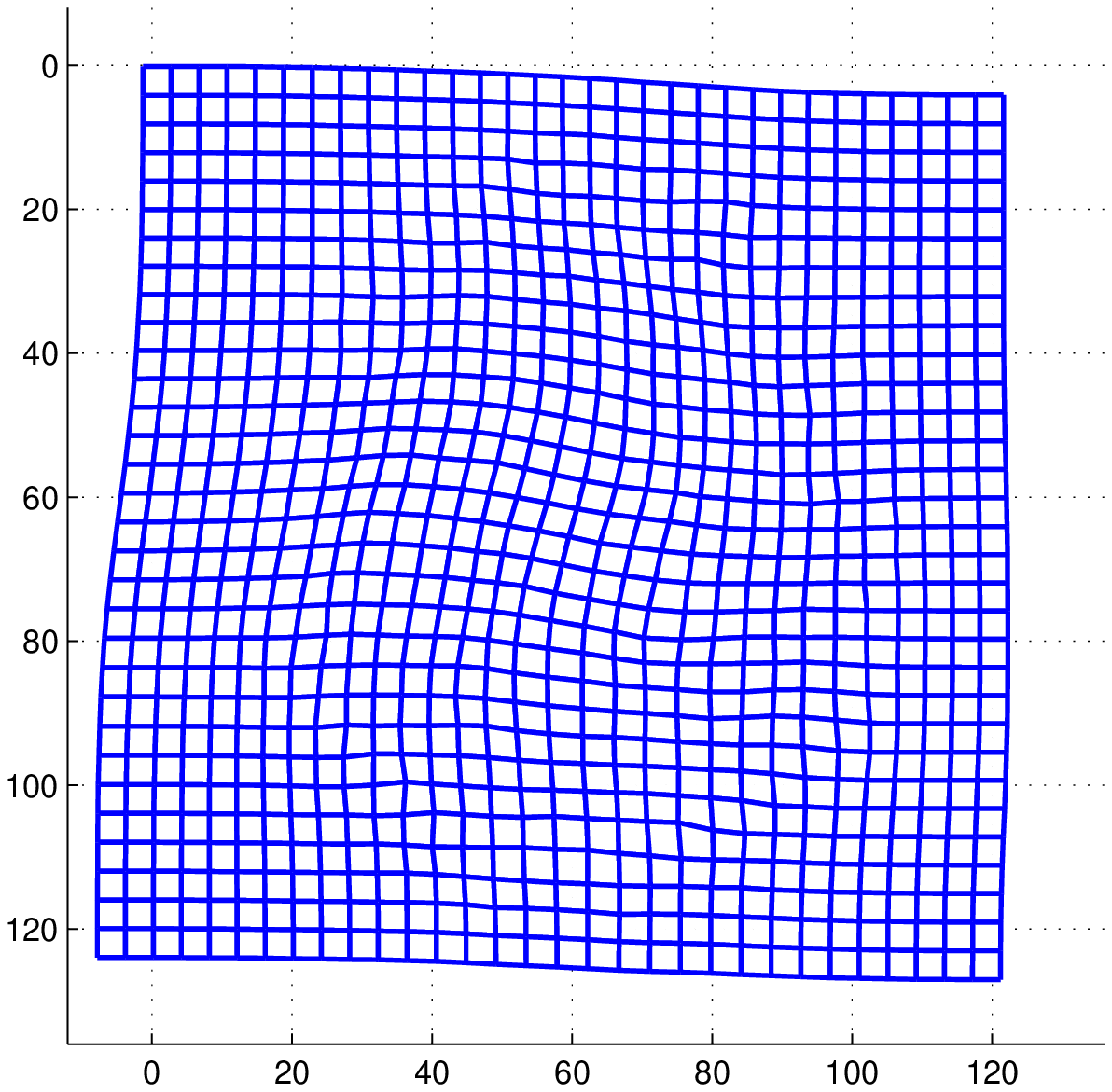}}
 \end{tabular}
\caption{Results of Gaussian curvature image registration for multi-modality images. The model is able to register multi-modality images with mutual information as the distance measure.}
\label{f1MM}
\end{figure}

{\bf A brief summary}.
The linear curvature model is relatively easy to solve, based on approximation of the mean curvature. The mean curvature model for image registration is highly nonlinear, making it challenging to solve. The Gaussian curvature resembles the mean curvature in many ways, though different,  but its model appears to deliver better quality
than the mean curvature.
 The diffeormorphic demon model is equivalent to the second order gradient descent on the SSD as shown in \cite{Pennec1999}. The model is only limited to mono-modality images and it is not yet applicable to multi-modality images. Our Gaussian curvature model however can be easily modified to work with  multi-modality images by replacing the SSD by a mutual information or normalised gradient fields based regularizer;
 an optimal solver is yet to be developed.
We show one example of extension to a pair of multi-modality images in Fig.\ref{f1MM}.

\section{Conclusions}
We have introduced a novel regularisation term for non-parametric image registration based on the Gaussian curvature of the surface induced by the displacement field. The model can be effectively solved using the augmented Lagrangian and alternating minimisation methods. For comparison, we used three models: the linear curvature \cite{Fischer2003}, the mean curvature
 \cite{CCC(2011)} and the demon algorithm \cite{Vercauteren2009} for   mono-modality images. Numerical experiments show that the proposed model delivers better results than the competing models.
%Numerical experiments show that the proposed model is more robust than these two models based on the tested images.
%\appendix
\section*{Appendix A -- Derivation of the Euler-Lagrange Equations} \label{AppA}
Let $q_1=u_x$ and $q_2=u_y$; then we can write the Gaussian curvature regularisation term as
\begin{equation}
{\mathcal S}^{GC}(q_1,q_2)=\int_{\Omega} \Big |\frac{q_{1,x} q_{2,y}-q_{1,y} q_{2,x} }{(1+q_1^2+q_2^2)^2}\Big| dx dy.
\nonumber
\end{equation}
From the optimality condition $\frac{d {\mathcal S}^{GC}(q_1,q_2)}{d q_1}=\frac{d {\mathcal S}^{GC}(q_1,q_2)}{d q_2}=0$, then $\frac{d}{d \epsilon_1}{\mathcal S}^{GC}(q_1+\epsilon_1 \varphi_1,q_2)\Big |_{\epsilon_1=0}=0$
and $ \frac{d}{d \epsilon_2}{\mathcal S}^{GC}(q_1,q_2+\epsilon_2 \varphi_2)\Big |_{\epsilon_2=0}=0$. In details,
\begin{equation}\label{EL1}
\mbox{ }\hspace*{-10mm}
\begin{aligned}
 \frac{d}{d \epsilon_1} \int_{\Omega} \Big|\frac{(q_{1}+\epsilon_1 \varphi_1)_x q_{2,y}-(q_{1}+\epsilon_1 \varphi_1)_y  q_{2,x} }{(1+(q_{1}+\epsilon_1 \varphi_1)^2+q_2^2)^2}\Big|  dx dy \Big |_{\epsilon=0}=\mbox{ }\hspace*{25mm}& \\
\int_{\Omega}  S \frac{d}{d \epsilon_1} \left[  \frac{(q_{1}+\epsilon_1 \varphi_1)_x q_{2,y}-(q_{1}+\epsilon_1 \varphi_1)_y  q_{2,x} }{(1+(q_{1}+\epsilon_1 \varphi_1)^2+q_2^2)^2}\right ]dx dy \Big |_{\epsilon=0}=0&
\end{aligned}
\end{equation}
where $S=\text{sign}\left( \frac{q_{1,x}q_{2,y}-q_{1,y}q_{2,x}}{(1+q_1^2+q_2^2)^2}  \right)$. From (\ref{EL1}),
\begin{equation}
\begin{aligned}
& \int_{\Omega} S \left[ \frac{\varphi_{1,x} q_{2,y} -\varphi_{1,y} q_{2,x}}{(1+q_1^2+q_2^2)^2}+(q_{1,x}q_{2,y}-q_{1,y}q_{2,x})(-4\varphi_1 q_1(1+q_1^2+q_2^2)^{-3}) \right]dx dy\\
&= \int_{\Omega} \frac{S\varphi_{1,x}q_{2,y}}{\Gamma^2}-\frac{S \varphi_{1,y}q_{2,x}}{\Gamma^2}-\frac{4SDq_1\varphi_1}{\Gamma^3}dx dy =0,\\
&\mbox{}\hspace{-1.2cm}\text{ where} \hspace{0.1cm} \Gamma=1+q_1^2+q_2^2, D=q_{1,x}q_{2,y}-q_{1,y}q_{2,x}.\\
\end{aligned}
\nonumber
\end{equation}
Using the Green theorem
$
\int_{\partial \Omega} \phi \boldsymbol{\omega} \cdot \boldsymbol{n} ds-\int_{\Omega} \phi \text{div} (\boldsymbol{\omega}) dx dy=\int_{\Omega} \nabla \phi \cdot \boldsymbol{\omega} dx dy,
$
we have,
\begin{equation}
\begin{aligned}
&\int_{\Omega} \frac{S\varphi_{1,x}q_{2,y}}{\Gamma^2}-\frac{S \varphi_{1,y}q_{2,x}}{\Gamma^2}dx dy=
\int_{\partial \Omega} \varphi_1 \left(\frac{Sq_{2,y}}{\Gamma^2}, \frac{Sq_{2,x}}{\Gamma^2}\right)  \cdot \boldsymbol{n} ds-\int_{\Omega} \varphi_1  \text{div}\left(\frac{Sq_{2,y}}{\Gamma^2}, \frac{Sq_{2,x}}{\Gamma^2}\right)=0
\end{aligned}
\nonumber
\end{equation}
 where $\phi=\varphi_1, \boldsymbol{\omega}=\left(\frac{Sq_{2,y}}{\Gamma^2}, \frac{Sq_{2,x}}{\Gamma^2}\right) $.
Setting the boundary integral to   zero, then we derive
\begin{equation}
\begin{aligned}
& \int_{\Omega} \varphi_1  \text{div}\left(\frac{Sq_{2,y}}{\Gamma^2}, \frac{Sq_{2,x}}{\Gamma^2}\right)=0.
\end{aligned}
\nonumber
\end{equation}
Finally, we use the fundamental lemma of calculus of variation to get:
\begin{equation}
\nabla \cdot \left(\frac{Sq_{2,y}}{\Gamma^2}, \frac{Sq_{2,x}}{\Gamma^2}\right)-\frac{4SDq_1}{\Gamma^3} =0.
\nonumber
\end{equation}
Similarly, for $\frac{d}{d \epsilon_2}{\mathcal S}^{GC}(q_1,q_2+\epsilon_2 \varphi_2)\Big |_{\epsilon_2=0}=0$, we finally obtain equation (\ref{EL}).
\hfill $\Box$

\bibliographystyle{abbrv}
\bibliography{Reference1_GC}
{\bf Authors' home page: \ \  \ \ http:$\backslash\backslash$www.liv.ac.uk/cmit}

\end{document}